\documentclass[reqno,11pt]{amsart}

\usepackage[left=1in,right=1in,top=1in,bottom=1in]{geometry}
\usepackage[T1]{fontenc}
\usepackage{lmodern}
\usepackage{microtype}
\microtypesetup{expansion=false}
\usepackage{enumitem}
\usepackage{amsmath}
\usepackage{amssymb}
\usepackage{mathtools}
\usepackage{aliascnt}
\usepackage{graphicx}
\usepackage{needspace}
\usepackage{placeins}
\usepackage{tikz}
\usetikzlibrary{arrows.meta}
\usepackage{array,booktabs}
\usepackage[colorlinks=true,citecolor=blue,linkcolor=blue,urlcolor=blue]{hyperref}
\usepackage[nameinlink,noabbrev]{cleveref}

\numberwithin{equation}{section}

\theoremstyle{plain}
\newtheorem{theorem}{Theorem}[section]
\newaliascnt{corollary}{theorem}

\aliascntresetthe{corollary}
\newaliascnt{con}{theorem}

\aliascntresetthe{con}
\newaliascnt{lemma}{theorem}
\newtheorem{lemma}[lemma]{Lemma}
\aliascntresetthe{lemma}
\newaliascnt{proposition}{theorem}
\newtheorem{proposition}[proposition]{Proposition}
\aliascntresetthe{proposition}
\newaliascnt{assum}{theorem}

\aliascntresetthe{assum}

\theoremstyle{definition}
\newaliascnt{definition}{theorem}
\newtheorem{definition}[definition]{Definition}
\aliascntresetthe{definition}
\newaliascnt{example}{theorem}

\aliascntresetthe{example}
\newaliascnt{q}{theorem}

\aliascntresetthe{q}

\theoremstyle{remark}
\newaliascnt{remark}{theorem}

\aliascntresetthe{remark}

\crefname{theorem}{theorem}{theorems}
\Crefname{theorem}{Theorem}{Theorems}
\crefname{corollary}{corollary}{corollaries}
\Crefname{corollary}{Corollary}{Corollaries}
\crefname{con}{conjecture}{conjectures}
\Crefname{con}{Conjecture}{Conjectures}
\crefname{lemma}{lemma}{lemmas}
\Crefname{lemma}{Lemma}{Lemmas}
\crefname{proposition}{proposition}{propositions}
\Crefname{proposition}{Proposition}{Propositions}
\crefname{assum}{assumption}{assumptions}
\Crefname{assum}{Assumption}{Assumptions}
\crefname{definition}{definition}{definitions}
\Crefname{definition}{Definition}{Definitions}
\crefname{example}{example}{examples}
\Crefname{example}{Example}{Examples}
\crefname{q}{question}{questions}
\Crefname{q}{Question}{Questions}
\crefname{remark}{remark}{remarks}
\Crefname{remark}{Remark}{Remarks}

\renewcommand{\theenumi}{\roman{enumi}}
\newlist{maincomments}{enumerate}{1}
\setlist[maincomments]{
  label=\textnormal{(\roman*)},
  ref=\roman*,
  leftmargin=*
}
\setlist[enumerate,1]{
  label=\textnormal{(\theenumi)},
  ref=\theenumi,
  leftmargin=*
}

\DeclarePairedDelimiterX{\inner}[2]{\langle}{\rangle}{#1,#2}

\newcommand{\R}{\mathbb R}

\newcommand{\dist}{\operatorname{dist}}

\newcommand{\Z}{\mathbb Z}
\newcommand{\Sph}{\mathbb S^2}
\newcommand{\HH}{\mathcal H}
\newcommand{\DD}{\mathcal D_1}
\newcommand{\axis}{\mathcal A}
\newcommand{\Sing}{\operatorname{Sing}}
\newcommand{\Reg}{\operatorname{Reg}}
\newcommand{\spt}{\operatorname{spt}}

\newcommand{\Mass}{\mathbf M}
\newcommand{\graph}[1]{\curr{\operatorname{graph}(#1)}}
\newcommand{\curr}[1]{[\![#1]\!]}

\hypersetup{
  pdftitle={Finiteness of singularities for higher-weight axially symmetric minimizers of the relaxed Dirichlet energy},
  pdfauthor={Zhiyuan Dai, Huaijie Wang, Wei Wang},
  pdfkeywords={harmonic maps; relaxed Dirichlet energy; minimal connections; Cartesian currents; axial symmetry; isolated singularities}
}

\title[Finiteness of singularities for axial minimizers]{Finiteness of singularities for higher-weight axially symmetric minimizers of the relaxed Dirichlet energy}
\author{Zhiyuan Dai}
\address{School of Mathematical Sciences, Zhejiang University, Hangzhou 310058, China}
\email{12635038@zju.edu.cn}
\email{zydai2026@outlook.com}

\author{Huaijie Wang}
\address{School of Mathematical Sciences, Peking University, Beijing 100871, China}
\email{huaijie\_wang@163.com}

\author{Wei Wang}
\address{School of Mathematical Sciences, Peking University, Beijing 100871, China}
\email{wwmath166@outlook.com}
\email{2201110024@stu.pku.edu.cn}
\date{}

\begin{document}

\begin{abstract}
Let $k\ge2$ be an integer. We prove that every minimizer of the
Bethuel--Brezis--Coron relaxed Dirichlet energy in the signed
meridional $k$-axial class, with continuous $H^1$ degree-zero data
on the boundary of a ball, is smooth outside a finite subset of the
symmetry axis. This extends the weight-one finiteness theorem of
Hardt--Lin--Poon and resolves the finite-singularity question posed
by Martinazzi--Rivi\`ere as Open Problem~A, in the fixed-boundary
Bethuel--Brezis--Coron setting with continuous degree-zero data.
The interior argument uses a distributional graph-boundary identity,
height-preserving half-turn replacements, and exact square
completions retaining the current cost. Cap-valued crosscuts and a
finite incidence count exclude singular accumulation at the boundary
poles.
\end{abstract}

\subjclass[2020]{Primary 58E20; Secondary 49Q20, 35J50}
\keywords{Harmonic maps, relaxed Dirichlet energy, minimal connections,
Cartesian currents, axial symmetry, isolated singularities}

\maketitle

\section{Introduction}\label{sec:intro}

The relaxed Dirichlet energy records both the energy of a Sobolev map
and the cost of the topological defects that can appear under weak
convergence. In three dimensions, this second contribution is measured
by the length of a connection between defects. The purpose of this paper is
to prove that, under higher-weight axial symmetry, its singularities
are finite in number when the prescribed boundary values
are continuous and have degree zero, extending the weight-one theorem of Hardt--Lin--Poon \cite{HLP}.

\subsection{From minimal connections to relaxed energy}

We use the normalization
\begin{equation}\label{eq:Dirichlet-energy}
 E(u;\Omega)=\int_\Omega|\nabla u|^2\,dx,
 \qquad u\in H^1(\Omega,\Sph),\quad\Omega\subset\R^3.
\end{equation}
The fixed-trace problem explains why a relaxation is needed. For smooth
boundary data $g$ of degree zero, set
\[
 \mathcal H_g=\{v\in H^1(B_R,\Sph):\operatorname{Tr}v=g\},\qquad
 \mathcal C_g=\mathcal H_g\cap C^0(\overline{B_R},\Sph).
\]
Although continuous extensions exist, minimization of $E$ in
$\mathcal H_g$ can produce singular maps and an energy gap relative
to $\mathcal C_g$ \cite{HardtLin}. Minimization directly in
$\mathcal C_g$ lacks weak compactness. The appropriate replacement is
\begin{equation}\label{eq:relaxation-motivation}
 \mathcal R_g(u)=\inf\left\{\liminf_{j\to\infty}E(u_j):
 u_j\in\mathcal C_g,\ u_j\rightharpoonup u\text{ in }H^1\right\}.
\end{equation}
The relaxation preserves the continuous-class infimum while admitting
weak limits:
\[
 \min_{\mathcal H_g}\mathcal R_g=\inf_{\mathcal C_g}E.
\]
Bethuel--Brezis--Coron identify $\mathcal R_g=E+8\pi L$
\cite{BBC}; see also Brezis's account \cite[Section~4]{BrezisProblems}.
Thus the connection term measures the energy required to approximate
a limiting map by continuous extensions. Lin \cite{LinLieb90}
surveys this connection between relaxed energies, defects, and minimal
currents. The approximation results of
Bethuel and Bethuel--Zheng \cite{Bethuel,BZ} place this phenomenon in
the broader density problem for manifold-valued Sobolev maps.

Brezis--Coron--Lieb \cite{BCL} identified the minimal-connection
structure associated with the distributional Jacobian
\[
 D(u)=\bigl(u\cdot(u_{x_2}\times u_{x_3}),
            u\cdot(u_{x_3}\times u_{x_1}),
            u\cdot(u_{x_1}\times u_{x_2})\bigr).
\]
It satisfies $|D(u)|\le |\nabla u|^2/2$; for finitely many point defects
of degrees $d_i$,
$\operatorname{div}D(u)=4\pi\sum_i d_i\delta_{a_i}$ in the interior.
The corresponding connection cost is $8\pi$ times its length.
Bethuel--Brezis--Coron \cite{BBC} incorporated this cost into the weak
relaxation of the fixed-boundary Dirichlet problem. For continuous
$H^1$ degree-zero data $g$ on $\partial B_R$, write
\begin{equation}\label{eq:L-definition}
 \begin{aligned}
 E_1(u;B_R)&=E(u;B_R)+8\pi L(u;B_R),\\
 L(u;B_R)&=\frac1{4\pi}
 \sup_{\substack{\xi\in C^1(\overline{B_R})\\|\nabla\xi|\le1}}
 \left\{\int_{B_R}D(u)\cdot \nabla\xi\,dx
       -\int_{\partial B_R}(\operatorname{Jac}g)\xi\,d\mathcal H^2\right\},
 \end{aligned}
\end{equation}
where $\operatorname{Tr}u=g$ and $\operatorname{Jac}g$ is the oriented tangential
Jacobian of $g$. Following \cite[Sections~4--5]{HLP}, $L(u)$ denotes the connection
length (minimal connection), $I_u$ its canonical mass-minimizing one-current, and
$I_T$ the vertical one-current of a general Cartesian current $T$.
In particular,
\[
 L(u)=\Mass(I_u),\qquad E_1(u)=E(u)+8\pi L(u).
\]
More generally, Hardt--Lin--Poon consider $E_\lambda=E+8\pi\lambda L$;
we study the critical case $\lambda=1$.
When the domain is fixed, we omit it from $E$, $E_1$, and $L$.

The liquid-crystal variational setting provides a related source of
existence and singularity questions \cite{HKL86,HKL}. For the
ordinary minimizing problem, the theory of
Schoen--Uhlenbeck \cite{SU,SU84} gives a locally finite singular set
in dimension three, with a separate boundary regularity theory
\cite{SU83}; degree-zero boundary values can still force
singular minimizers \cite{HardtLin}. Zhang \cite{Zhang1989} obtained
regular harmonic extensions for a class of nonsurjective axial data
and exhibited extensions that are not Dirichlet minimizers. This
distinguishes smooth harmonic extension from energy minimization. For the relaxed problem, a
competitor can trade graph energy against connection length. At a
subcritical coefficient, strict separation of those costs supplies
additional compactness \cite{BB,HLP}. At the critical coefficient
in \eqref{eq:L-definition}, that separation is lost; see also
\cite{Hardt}. The Cartesian-current formulation of
Giaquinta--Modica--Sou\v{c}ek \cite{GMS,GMSDirichlet,GMSBook}
uses the language of geometric measure theory \cite{Federer} and
makes the graph and connection energies local in a single object,
$T=\graph{u}+I_T\times\curr{\Sph}$. We use this formulation to construct
and compare competitors.

The dipole constructions originate in the work of Brezis--Coron
\cite{BC} and Brezis--Coron--Lieb \cite{BCL}.
Hardt--Lin--Poon \cite{HLP} solved the finite-singularity problem in
the weight-one meridional axial class. In particular, their
Theorem~9.1 gives global finiteness for continuous $H^1$ degree-zero
boundary values at the critical coefficient. They also prove
$L(u)=0$ and degree zero at every isolated interior singularity.
Their argument combines a strict dipole comparison, the geometry of
pole-level curves, and smooth replacements; the main comparison
results are \cite[Lemma~7.1 and Theorem~8.2]{HLP}.

Two difficulties arise at higher axial weight. First, a smooth
weight-$k$ map has transverse component $O(r^k)$ near the symmetry
axis. When $k\ge2$, its derivative vanishes there, so the hypothesis
$\nabla u(0)\ne0$ in Hardt--Lin--Poon's strict dipole lemma is unavailable. Second,
Martinazzi \cite[Theorems~1 and~2]{Martinazzi} constructed a smooth,
nonconstant higher-weight graph minimizing with a nonzero vertical
current in the corresponding structural class. Vertical mass thus
cannot be treated as a graph singularity or excluded solely by
regularity. Martinazzi--Rivi\`ere \cite[p.~1222, Open Problem~A]{MR}
ask whether higher-weight axially symmetric minimizers of the relaxed
energy are smooth away from finitely many points. Theorem~\ref{thm:global}
answers this question affirmatively and resolves Open Problem~A in
the fixed-boundary Bethuel--Brezis--Coron setting with continuous $H^1$ degree-zero data.
The precise boundary assumptions are included in the theorem:
its conclusion is global finiteness of the interior singular set.
The proof first establishes interior local finiteness for structural
currents with compatible vertical baselines in $k\Z$.

These are symmetry-restricted statements. Brezis's Problems~4.2,
4.3, and~4.4 \cite[Section~4]{BrezisProblems} ask, respectively,
about attainment among continuous extensions, continuity of relaxed
minimizers, and vanishing of the connection term in the unrestricted
class. The results here do not identify axial minimizers with
unrestricted minimizers, and the global theorem asserts finiteness
rather than vanishing of $L(u)$.

\subsection{Main theorem}

Write $x=(r\cos\theta,r\sin\theta,z)$ and
$\axis=\{r=0\}$, oriented by increasing $z$.
For an integer $k\ge1$, a signed meridional $k$-axial map is
represented on each simply connected off-axis meridian patch by
\begin{equation}\label{eq:ansatz}
 u(r,\theta,z)=
 \bigl(\cos(k\theta)\sin\varphi(r,z),\,
       \sin(k\theta)\sin\varphi(r,z),\,\cos\varphi(r,z)\bigr).
\end{equation}
The real-valued function $\varphi$ is called an \emph{angle function}
for $u$, as in \cite[Section~1]{HLP}. No sign restriction is imposed
on $\sin\varphi$. Angle functions for the same circle profile differ
by constant multiples of $2\pi$ on connected overlaps. The functions
$\sin\varphi$ and $\cos\varphi$ are therefore globally defined, even
when a single real angle function on the whole meridian is not.
On a ball its open meridian is simply connected and one may choose
a global angle function there; the Sobolev angle construction is given in
Subsection~\ref{sec:scalar-equation}.
For boundary data $g$ of this form, set
\[
 \mathcal A_k(g)=\{u\in H^1(B_R,\Sph):
       u\text{ has \eqref{eq:ansatz}},\ \operatorname{Tr}u=g\}.
\]
Let $\Reg(u)$ denote the open set on which the Sobolev class of $u$
has a smooth representative, and put
$\Sing(u)=\Omega\setminus\Reg(u)$ for a map defined on $\Omega$.

\begin{theorem}[Global finiteness with continuous boundary values]
\label{thm:global}
Let $k\ge2$, $R>0$, and let
$g\in H^1(\partial B_R,\Sph)\cap C^0(\partial B_R,\Sph)$ be
signed meridional $k$-axial with $\deg g=0$.
The functional $E_1$ in \eqref{eq:L-definition} attains its minimum
in $\mathcal A_k(g)$. Every minimizer satisfies
\[
 \Sing(u)\subset\axis\cap B_R,\qquad \#\Sing(u)<\infty,
 \qquad u\in C^\infty(B_R\setminus\Sing(u),\Sph).
\]
\end{theorem}

This extends the finiteness assertion of
\cite[Theorem~9.1(2)]{HLP} to every higher axial weight. Global
finiteness requires control of possible accumulation at the two
boundary poles; the proof supplies that control without asserting
regularity at those poles.

\subsection{Proof strategy and relation to earlier methods}

The proof retains Hardt--Lin--Poon's local current formulation, pole-level
geometry in the meridian half-plane, and analytic continuation
after an energy-preserving replacement. The changes needed at
higher weight enter at the comparison steps.

We first determine the graph boundary before making any
bounded-variation assumption on the axis trace:
\[
 \partial \graph{u}=-\frac{k}{2}q'\times\curr{\Sph},\qquad q+2m=c.
\]
Subtracting a constant vertical baseline reduces the line cost to
two adjacent states. A fold into a small polar cap then gives an
open dense set of regular axis points. The vertical current is
kept as a separate part of the comparison throughout.

For the pole-level topology, we replace $u$ inside a meridian Jordan
domain by $R_\pi u$, where
$R_\pi=\operatorname{diag}(-1,-1,1)$, and leave $I_T$ fixed.
This half-turn preserves the height, graph energy, and current cost.
Equality and analytic continuation exclude pole-level cycles and
arcs with two interior axis ends. A finite-tree argument then
controls all axis contacts, including limits of whole graph
components, and proves their accessibility by nodal arcs.

On each remaining one-sided angle region, two square completions
give comparisons with
\[
 2\arctan\left(\frac{\varepsilon}{r^k}\right),\qquad
 2\arctan\left(\frac{r^k}{\varepsilon}\right).
\]
The first changes the axis trace and the current together. The
second is used after the axis trace has been fixed. Their equality
cases imply the order-$r^k$ estimate required for smooth removal.
These comparisons give at most two singularities on a one-sided
axis interval. Together with the finite pole-contact count, they
prove Lemma~\ref{lem:interior}.

For Theorem~\ref{thm:global}, continuity of the prescribed trace
provides a polar-cap collar near each physical boundary pole.
Following Hardt--Lin--Poon's boundary localization, we place an artificial
interior segment across that collar. A cap-valued crosscut permits
a trace-preserving fold that clears an axis tail. If no such
crosscut exists, the half-turn comparison forces every remaining
pole-level contact to reach a distinct incidence on the artificial
segment. There are finitely many incidences. The one-sided
comparison therefore bounds the singularities near each boundary
pole, completing the global count.

Section~\ref{sec:scalar} establishes the local identities and cap
comparisons. Section~\ref{sec:topology} proves interior local finiteness
by combining the finite contact count with the two-point comparison.
Section~\ref{sec:map} identifies the canonical connection, proves
existence, and transfers minimality to currents.
Section~\ref{sec:boundary} controls the two boundary poles and completes
the proof of Theorem~\ref{thm:global}.

\section{Scalar identities and cap comparisons}\label{sec:scalar}

We work with structural Cartesian currents, following the local formulation
of Hardt--Lin--Poon \cite[Section~5]{HLP}.
For a rotationally invariant open set $\Omega\subset\R^3$, let
$\operatorname{cart}_{\mathrm{AS}}^{(k)}(\Omega,\Sph)$ consist of
\begin{equation}\label{eq:cartesian}
 \begin{aligned}
 T&=\graph{u}+I_T\times\curr{\Sph},\qquad
 I_T=km(z)\curr{\axis},\qquad
 m\in L^1_{\rm loc}(\axis\cap\Omega,\Z),\\
 &u\in H^1_{\rm loc}(\Omega,\Sph)\text{ of the form \eqref{eq:ansatz}},
 \qquad \partial T=0\text{ in }\Omega\times\Sph,
 \end{aligned}
\end{equation}
with locally finite energy
\begin{equation}\label{eq:energy}
 \DD(T;U)=E(u;U)+8\pi\Mass(I_T \llcorner U),\qquad U\Subset\Omega.
\end{equation}
Here $\graph{u}$ is the oriented Sobolev graph current. Its
first and second minors are locally $L^2$ and $L^1$; its third minors
vanish because the target is $\Sph$. The graph has multiplicity-one
horizontal projection and no mass over the axis. These conditions
follow the Cartesian-current framework of \cite{GMS} and the axial
restriction in \cite{Martinazzi}.
The superscript $(k)$ records the higher axial weight. The class here
also allows constant vertical baselines in $k\Z$; the normalization
in Subsection~\ref{sec:normalization} reduces them to two adjacent cost states.

\begin{definition}\label{def:local-minimality}
A current $T\in\operatorname{cart}_{\mathrm{AS}}^{(k)}(\Omega,\Sph)$ is locally minimizing if
\[
 \DD(T;U)\le\DD(\widetilde T;U)
\]
for every rotationally invariant open $U\Subset\Omega$ and every
$\widetilde T\in\operatorname{cart}_{\mathrm{AS}}^{(k)}(\Omega,\Sph)$ with
$\spt(\widetilde T-T)\Subset U\times\Sph$.
Both the graph and the vertical current may change.
\end{definition}

Throughout the interior argument, $T$ is locally minimizing.
The axis trace is initially only measurable; no bounded-variation
assumption enters the comparisons.

\subsection{The scalar equation and off-axis regularity}\label{sec:scalar-equation}

For an open set $V$ in the meridian half-plane $\{(r,z):r>0\}$,
let $\widehat V\subset\R^3$ be its solid of revolution, with the
interior axis points included:
\[
 \begin{split}
 \widehat V={}&\{(r\cos\theta,r\sin\theta,z):
       (r,z)\in V,\ 0\le\theta<2\pi\}\\
 &\cup\{(0,0,z):\exists\varepsilon>0,
 (0,\varepsilon)\times(z-\varepsilon,z+\varepsilon)\subset V\}.
 \end{split}
\]
In particular, the meridian rectangle $C=(0,R)\times I$ corresponds
to the open cylinder $\widehat C=B_R^2\times I$, where
$B_R^2=\{(x_1,x_2):x_1^2+x_2^2<R^2\}$.
Energies on meridian regions are evaluated on their solids of
revolution, including the axis-current contribution.

We will use the following elementary extension argument. If a bounded
function $v$ is locally $H^1$ off the axis and its off-axis gradient
$g$ lies in $L^2(U)$ for a bounded open set $U$, then $v\in H^1(U)$
with $\nabla v=g$. Indeed, set
\[
 \eta_\varepsilon(r)=
 \begin{cases}
 0,&r\le\varepsilon^2,\\
 \log(r/\varepsilon^2)/|\log\varepsilon|,
       &\varepsilon^2<r<\varepsilon,\\
 1,&r\ge\varepsilon.
 \end{cases}
 \qquad
 \int_U|\nabla\eta_\varepsilon|^2\le C_U|\log\varepsilon|^{-1}.
\]
Testing off the axis with $\eta_\varepsilon\zeta$, for
$\zeta\in C_c^\infty(U)$, gives
\[
 \int\eta_\varepsilon v\,\partial_i\zeta
 =-\int\eta_\varepsilon g_i\zeta
  -\int v\zeta\,\partial_i\eta_\varepsilon,
 \qquad
 \left|\int v\zeta\,\partial_i\eta_\varepsilon\right|
 \le\|v\|_\infty\|\zeta\|_2\|\nabla\eta_\varepsilon\|_2\to0.
\]
Passing to the limit identifies the weak gradient on all of $U$.

The meridional representation determines an angle modulo $2\pi$;
we first justify an $H^1$ choice. On an off-axis meridian disk $D$,
angular averaging recovers
\[
 w=\frac1{2\pi}\int_0^{2\pi}
 (u_1\cos k\theta+u_2\sin k\theta,u_3)\,d\theta
 \in H^1(D;S^1).
\]
For $j=w_2\,dw_1-w_1\,dw_2$, the circle constraint gives
\[
 dj=2\,dw_2\wedge dw_1=0,\qquad
 dw_1=w_2j,\quad dw_2=-w_1j.
\]
These identities hold distributionally by smooth approximation.
The weak Poincare lemma yields $j=d\varphi$; then
$d(e^{-i\varphi}(w_2+iw_1))=0$, so an additive adjustment gives
$w=(\sin\varphi,\cos\varphi)$ with $\varphi\in H^1(D)$.
On a simply connected meridian the same argument gives a global
$H^1_{\rm loc}$ angle function.

Cylindrical integration gives
\begin{equation}\label{eq:scalar-energy}
 E(u;\widehat C)=2\pi\int_C
 \left[r(\varphi_r^2+\varphi_z^2)
       +\frac{k^2\sin^2\varphi}{r}\right]dr\,dz.
\end{equation}
For $\eta\in C_c^\infty(D)$, the variation
$\varphi_t=\varphi+t\eta$ is induced by smooth fiber rotations
supported away from the axis. It preserves closure and leaves $I_T$
unchanged. Local minimality therefore gives
\[
 0=\frac1{4\pi}\left.\frac d{dt}E(u_t)\right|_{t=0}
 =\int_D\left(r\nabla\varphi\cdot\nabla\eta
       +\frac{k^2}{r}\sin\varphi\cos\varphi\,\eta\right)dr\,dz.
\]
Thus
\begin{equation}\label{eq:pde}
 \varphi_{rr}+r^{-1}\varphi_r+\varphi_{zz}
       -k^2r^{-2}\sin\varphi\cos\varphi=0.
\end{equation}
Interior elliptic estimates and bootstrapping, followed by analytic
regularity \cite{Morrey}, give $\varphi,u\in C^\omega$ off the axis.
Direct substitution into \eqref{eq:ansatz} gives the vector equation
there. To extend it across the axis, test with
$\eta_\varepsilon\zeta$, where
$\zeta\in C_c^\infty(\Omega;\R^3)$. The only additional term satisfies
\[
 \left|\int\nabla u:(\zeta\otimes\nabla\eta_\varepsilon)\right|
 \le\|\zeta\|_\infty
 \|\nabla u\|_{L^2(\spt\zeta)}
 \|\nabla\eta_\varepsilon\|_{L^2(\spt\zeta)}\to0.
\]
Dominated convergence yields
\begin{equation}\label{eq:weak-harmonic}
 -\Delta u=|\nabla u|^2u\quad\text{in }\mathcal D'(\Omega;\R^3).
\end{equation}
Uniformly bounded smooth approximation also permits compactly
supported $H^1\cap L^\infty$ test functions.

\subsection{The axis trace and graph-boundary identity}
\label{subsec:axis-boundary}

Finite energy supplies a pole-valued axis trace. Indeed, on
$C=(0,R)\times I$, for $\ell=r,z$,
\[
 \int_C|\partial_\ell\cos\varphi|
 \le\left(\int_Cr\varphi_\ell^2\right)^{1/2}
     \left(\int_C\frac{\sin^2\varphi}{r}\right)^{1/2}<\infty.
\]
Hence $\cos\varphi\in W^{1,1}(C)$ has axis trace $q$.
For almost every $z$, the radial slice has an integrable derivative
and therefore a limit. If that limit lay in $(-1,1)$, its angular
energy $\int_0^R\sin^2\varphi\,dr/r$ would be infinite. Thus
\begin{equation}\label{eq:axis-trace}
 \cos\varphi(r,z)\to q(z)\in\{-1,1\},\qquad
 \sin\varphi(r,z)\to0\quad\text{for a.e. }z.
\end{equation}
These limits agree with the Sobolev trace. No trace of $\varphi$
itself, or bounded variation of $q$, is assumed.

\begin{lemma}[Graph-boundary identity]\label{lem:full-boundary}
For every finite-energy meridional $k$-axial map on
$\widehat C=B_R^2\times I$,
\begin{equation}\label{eq:full-boundary}
 \partial\graph{u}=-\frac{k}{2}q'\times\curr{\Sph}
 \quad\text{in }\widehat C\times\Sph.
\end{equation}
Consequently, for $I_T=km\curr{\axis}$ with integer-valued
$m\in L^1_{\rm loc}(I)$, the current
$T=\graph{u}+I_T\times\curr{\Sph}$ is closed precisely when
\begin{equation}\label{eq:c}
 q+2m=c\quad\text{a.e. on }I
\end{equation}
for an odd integer $c$. Derivatives are understood distributionally.
\end{lemma}
\begin{proof}
Smooth $H^1$ approximation of local angle functions gives convergence
of graph minors of orders one and two in $L^2$ and $L^1$; those of
order three vanish. Thus $\partial\graph{u}=0$ off the axis.

For a compactly supported smooth two-form $\eta$, choose
$\chi_\epsilon=0$ for $r\le\epsilon$ and $\chi_\epsilon=1$ for
$r\ge2\epsilon$, with $0\le\chi_\epsilon'\le C/\epsilon$.
Expanding $\graph{u}(d(\chi_\epsilon\eta))=0$ gives
\begin{equation}\label{eq:averaged-tube-boundary}
 \partial\graph{u}(\eta)
 =\lim_{\epsilon\downarrow0}\int\chi_\epsilon'(r)K_\eta(r)\,dr,
 \quad
 K_\eta(r)=\int_I\int_0^{2\pi}
 \eta_{(x,u)}((e_z,u_z),(re_\theta,u_\theta))\,d\theta\,dz.
\end{equation}
The tangent order is the inner tube orientation.
Write the vertical part of $\eta$ as $b(x,y)\omega(y)$, with
$\omega$ the oriented sphere area form, and freeze $x$ at $(0,0,z)$.
The horizontal, mixed, and freezing errors have cutoff average at most
\[
 C\epsilon+C\left(
 \int_\epsilon^{2\epsilon}\int_I
 \left[r\varphi_z^2+\frac{\sin^2\varphi}{r}\right]dz\,dr
 \right)^{1/2}\longrightarrow0.
\]
This follows from $|u_z|=|\varphi_z|$,
$|u_\theta|=k|\sin\varphi|$, and Cauchy--Schwarz.

For the surviving vertical term, put
\[
 \begin{split}
 A(z,s)&=\int_0^{2\pi}b\bigl((0,0,z),
 (\sqrt{1-s^2}\cos\gamma,\sqrt{1-s^2}\sin\gamma,s)\bigr)d\gamma,\\
 \mathcal B(z,s)&=k\int_s^1A(z,t)\,dt,\qquad
 W(z)=\int_{\Sph}b((0,0,z),y)\omega(y).
 \end{split}
\]
Since $u\cdot(u_z\times u_\theta)=-k\partial_z\cos\varphi$,
angular averaging and the Sobolev chain rule give
\[
 K_{\rm vert}(r)
 =-k\int_I A(z,\cos\varphi)\partial_z\cos\varphi\,dz
 =-\int_I\mathcal B_z(z,\cos\varphi)\,dz.
\]
The total derivative vanishes by compact support. Trace convergence
in $L^1(I)$ and $\mathcal B(z,1)=0$, $\mathcal B(z,-1)=kW(z)$ yield
\[
 \partial\graph{u}(\eta)
 =-\frac{k}{2}\int_I(1-q)W'\,dz
 =-\frac{k}{2}\langle q',W\rangle.
\]
Finally, $\partial I_T=-km'$, so closure is equivalent to
$(q+2m)'=0$. Since $q=\pm1$ and $m\in\Z$, its constant is odd.
\end{proof}

The Jacobian vector used in \eqref{eq:L-definition} is
\begin{equation}\label{eq:jacobian-vector}
 D(u)=\frac{k}{r}
 \bigl(\partial_z(\cos\varphi)e_r-\partial_r(\cos\varphi)e_z\bigr),
.
\end{equation}
Identity~\eqref{eq:full-boundary} applies to every finite-energy
meridional competitor, including the folds and clamps below.

The next observation makes nonstrict comparisons useful: an
energy-preserving replacement is again minimizing, so agreement on
an exterior open set propagates by analyticity.

\begin{lemma}[Equality in a compact comparison]\label{lem:equality}
Let $T$ be locally minimizing and $\widetilde T$ admissible, with
\[
 \spt(\widetilde T-T)\Subset U\times\Sph,\qquad
 \DD(\widetilde T;U)\le\DD(T;U),\qquad U\Subset\Omega.
\]
Then equality holds and $\widetilde T$ is locally minimizing.
Their graph maps coincide on every connected off-axis component
on which they agree on a nonempty open set. In particular, they
coincide in $B\setminus\axis$ if the change is compact in an
axis-centered ball $B\Subset\Omega$.
\end{lemma}
\begin{proof}
Minimality gives equality. For any compact perturbation $S$ of
$\widetilde T$, take a rotationally invariant $W\Subset\Omega$
containing $\overline U$ and its support. Additivity gives
\[
 \DD(\widetilde T;W)=\DD(T;W)\le\DD(S;W).
\]
Cancelling the common exterior energy proves local minimality.
Both graph maps are analytic off the axis; the identity theorem
proves the remaining assertions.
\end{proof}

\subsection{Current normalization and small-cap removal}\label{sec:normalization}

Normalization means subtracting the constant vertical current
common to the two possible axis phases. Work in a connected
interior axis neighborhood, where Lemma~\ref{lem:full-boundary} gives
$q+2m=c$ for an odd integer $c$. Set
\[
 P=\operatorname{sign}c,\qquad
 m_0=\frac{c-P}{2},\qquad
 T^0=T-km_0\curr{\axis}\times\curr{\Sph}.
\]
The graph map is unchanged, while the new vertical multiplicity is
\[
 m^0=m-m_0=\frac{P-q}{2}\in\{0,P\},
 \qquad q+2m^0=P.
\]
Thus the phase $q=P$ has zero normalized current, and the phase
$q=-P$ has multiplicity $Pk$.

The subtracted current has zero interior boundary. Moreover,
$|m|=|m_0|+|m^0|$, so on a comparison cylinder
$\widehat C=B_R^2\times J$,
\[
 \DD(T;\widehat C)
 =\DD(T^0;\widehat C)+8\pi k|m_0|\mathcal H^1(J).
\]
Every compactly supported competitor has the same constant $c$,
hence the same background $km_0\curr{\axis}$ and the same additive
energy term. Normalization therefore preserves closure and local
minimality.

Hereafter we omit the superscript $0$ and use the normalized notation
\begin{equation}\label{eq:normalized-measure}
 c=P\in\{-1,1\},\qquad m\in\{0,P\},\qquad
 \mu=|\nabla u|^2\,dx+8\pi k\HH^1\llcorner I,
\end{equation}
where $I=\{m\ne0\}$ is represented by its density-one points.
The original vertical current is recovered by adding
$km_0\curr{\axis}$.

\begin{lemma}[Cap-valued removal]\label{lem:cap-removal}
Let $k\ge2$. Suppose a finite-energy meridional map is smooth and harmonic
off the axis in a ball centered on the axis. If its off-axis image lies
in a polar cap of angular radius at most $\pi/4$, it is smooth on a
neighborhood of the center.
\end{lemma}
\begin{proof}
After a target reflection, choose the signed angle $|f|\le\alpha=\pi/4$.
The extension argument and finite energy give $f\in H^1$ physically
and $\int f^2/r<\infty$. Off the axis,
\[
 \mathcal Lf=0,\qquad
 \mathcal L=\partial_{rr}+r^{-1}\partial_r+\partial_{zz}
 -\frac{k^2a}{r^2},\qquad
 a=\frac{\sin f\cos f}{f}\ge\frac2\pi,
\]
with $a=1$ at $f=0$.
Translate the center to $z=0$ and choose a spindle compactly inside
the ball:
\[
 d(z)=\frac{H^2-z^2}{H},\quad \kappa=\frac14,\quad
 D=\{|z|<H,\ 0<r<\kappa d(z)\},\quad
 B=\frac{\alpha r}{\kappa d(z)}.
\]
On its curved side $B=\alpha$, while direct differentiation gives
\[
 \frac{(d^{-1})''}{d^{-1}}\le\frac{10}{d^2},\qquad
 \frac{\mathcal LB}{B}
 \le\frac{1-2k^2/\pi+10\kappa^2}{r^2}<0.
\]
Both $f$ and $B$ have finite weighted gradient and potential norms;
for $B$, radial integration reduces them to bounded multiples of
$1$ and $|d'|^2$, including at the tips.
The bounded function $w=(f-B)_+$ has zero side trace; the two tips
have zero $H^1$ capacity. Test the difference inequality with
$\eta_\varepsilon^2w$. The logarithmic cutoff error is bounded by
\[
 C\left(\int_{D\cap\{r<\varepsilon\}}
 r|\nabla(f-B)|^2\right)^{1/2}
 \left(\int_Drw^2|\nabla\eta_\varepsilon|^2\right)^{1/2}\to0.
\]
Thus
\[
 \int_D\left(r|\nabla w|^2+\frac{k^2a}{r}w^2\right)\le0.
\]
Applying the same argument to $-f$ gives $|f|\le B$, hence
$|f|\le Cr$ on a smaller central cylinder. For the physical scalar
$\Phi(x)=f(r,z)$, interior Poisson estimates on balls of radius
comparable to $r$ yield
\[
 |\nabla\Phi|
 \le C\bigl(r^{-1}\|\Phi\|_\infty+r\|k^2a\Phi/r^2\|_\infty\bigr)
 \le C,\qquad k|\sin f|/r\le C.
\]
Consequently $\nabla u$ is bounded. Extend the weak harmonic-map
equation across the axis by the logarithmic cutoff; its bounded
right side gives $W^{2,p}$ regularity and smoothness by bootstrapping.
\end{proof}

\subsection{Cap folds and dense regularity on the axis}\label{sec:cap-folds}

Fix an axis center $a=(0,0,a_z)$. For almost every sphere radius $s$,
the tangential energy is
\begin{equation}\label{eq:sphere-energy}
 S(s)=\int_{\partial B_s(a)}|\nabla_{\mathrm{tan}} u|^2
 =2\pi\int_0^\pi
 \left(\sin t\,\varphi_t^2+\frac{k^2\sin^2\varphi}{\sin t}\right)dt
 \ge4\pi k\int_0^\pi|\partial_t(\cos\varphi)|\,dt.
\end{equation}
The endpoint values of $\cos\varphi$ on this sphere are $q(a_z+s)$ and
$q(a_z-s)$ for almost every $s$. Here no axis regularity is assumed.
Indeed, in any fixed annulus $0<s_0<s<s_1$, the change of variables
\[
 (s,t)\longmapsto(r,z)=(s\sin t,a_z+s\cos t)
\]
has Jacobian of absolute value $s$, including at $t=0$ and $t=\pi$.
Near either endpoint it is a smooth coordinate map of manifolds with
boundary. The $W^{1,1}$ trace theorem and slicing applied to
$\cos\varphi\in W^{1,1}$ identify the endpoint limits with the axis trace.
Formula~\eqref{eq:sphere-energy} also gives absolute continuity of
$\cos\varphi$ along the semicircle. Since its endpoint values are poles,
$\sin\varphi\to0$ there, so the full sphere trace is continuous at
the two poles as well as smooth elsewhere.

If one endpoint has pole height $\sigma$, then throughout that sphere
\begin{equation}\label{eq:sphere-cap}
 1-\sigma \cos\varphi\le\frac{S(s)}{4\pi k}.
\end{equation}

For a polar cap of radius $0<\alpha<\pi/2$ about a pole $\sigma$, define
$\Pi_{\sigma,\alpha}:\Sph\to\Sph$ by preserving longitude and replacing
the colatitude $\vartheta$ from $\sigma$ by
\[
 \tau(\vartheta)=\min\{\alpha,\vartheta,\pi-\vartheta\}.
\]
This map fixes the cap and sends both poles to $\sigma$.
Since $|\tau'|\le1$ and $\sin\tau\le\sin\vartheta$,
it is one-Lipschitz and does not increase Dirichlet energy.

\begin{lemma}[The admissible cap fold]\label{lem:cap-fold}
Suppose a good sphere $\partial B_s(a)$ has image in a cap of
radius $\alpha_0<\alpha\le\pi/4$ about $\sigma$. Assume either
\[
 \text{(i)}\ q=\sigma\text{ a.e. on its axis interval},\qquad
 \text{or (ii)}\ c=P=\sigma\text{ after baseline normalization}.
\]
Then $u$ is smooth throughout $B_s(a)$. In case (ii), the normalized vertical current
vanishes on that interval.
\end{lemma}
\begin{proof}
Define $v=\Pi_{\sigma,\alpha}\circ u$ in $B_s(a)$ and $v=u$
outside. The fold fixes the sphere trace; Lipschitz composition
and Sobolev pasting give $v\in H^1$. On almost every radial
slice its axis limit is $\sigma$, so $q_v=\sigma$ on the interval.
Choose
\[
 m_v=\begin{cases}m,&\text{in case (i)},\\0,&\text{in case (ii)}
 \end{cases}
\]
on that interval and keep $m_v=m$ outside. In both cases
$q_v+2m_v=c$ globally. Lemma~\ref{lem:full-boundary} proves closure,
including the two sphere-axis points. With $T_v=\graph{v}+km_v
\curr{\axis}\times\curr{\Sph}$, one has
\[
 0\le\DD(T;B_s)-\DD(T_v;B_s)
 =E(u;B_s)-E(v;B_s)+8\pi k\int(|m|-|m_v|),
\]
where both summands on the right are nonnegative. The support of
the change is compact in a slightly larger axis-centered ball.
Minimality gives the reverse inequality. Therefore
\[
 E(u;B_s)=E(v;B_s),\qquad
 \int(|m|-|m_v|)=0,\qquad u=v\text{ off the axis}
\]
by Lemma~\ref{lem:equality}. Since $v(B_s)\subset
\overline{B_\alpha(\sigma)}$, the same inclusion holds for $u$
off the axis. Lemma~\ref{lem:cap-removal} proves regularity at every axis point inside the ball.
In case (ii), $m_v=0$ and the last integral gives $m=0$ a.e.
\end{proof}

\begin{lemma}[Small-energy regularity]\label{lem:small-energy}
There is $\epsilon_k>0$ with the following properties for every
axis-centered ball compactly contained in the minimizing region.
\begin{enumerate}
\item If $q$ and $m$ are constant on its axis interval and
$E(B_R(a))/R\le\epsilon_k$, then the center is regular.
\item After normalization, if $\mu(B_R(a))/R\le\epsilon_k$, then
the center is regular and the normalized current vanishes on an axis
neighborhood of it.
\end{enumerate}
\end{lemma}
\begin{proof}
Set $\alpha_0=\pi/8$, $\alpha=\pi/4$,
$\delta_0=1-\cos(\pi/8)$, and
$\epsilon_k=(\pi k/2)\delta_0$. Coarea gives
\[
 \int_{R/2}^R S(s)\,ds\le E(B_R(a)).
\]
In the first case let $\sigma$ be the constant trace pole; it may
be the high pole $-P$. Almost every radius has both endpoint traces
equal to $\sigma$. There is a good radius with $S(s)\le4\epsilon_k$.
By \eqref{eq:sphere-cap},
$1-\sigma \cos\varphi\le\delta_0/2<\delta_0$, so its image lies strictly in the
$\alpha_0$ cap. Fold and retain the constant current.

In the second case, the radii $s\in(R/2,R)$ for which either
$a_z+s$ or $a_z-s$ belongs to the high phase have measure at most
\[
 \HH^1(I\cap(a_z-R,a_z+R))
 \le\frac{\epsilon_kR}{8\pi k}
 =\frac{\delta_0R}{16}<\frac R4.
\]
Thus the good-pole radii have measure greater than $R/4$. Coarea
again supplies one with $S(s)\le4\epsilon_k$, after discarding all
slicing null sets. Its two poles are $\sigma=P=c$, and the same cap bound holds. Fold and set the normalized current to zero.
Lemma~\ref{lem:cap-fold} proves both assertions.
\end{proof}

\begin{lemma}[Dense regular axis]\label{lem:dense}
The regular points form an open dense subset of the interior axis.
\end{lemma}
\begin{proof}
Work with the normalized current, and let $I=\{m\ne0\}$ be
represented by its density-one points. For every compact interior
axis interval $J$, Fubini gives
\begin{equation}\label{eq:graph-fubini}
 \int_J\frac{E(B_\rho((0,0,z)))}{\rho}\,dz
 \le2\int_{\{r<\rho,\ \dist(x_3,J)<\rho\}}|\nabla u|^2\,dx
 \longrightarrow0.
\end{equation}
The limit follows from integrability and the shrinking tube.
Fatou's lemma along a sequence $\rho\downarrow0$ yields
\[
 \liminf_{\rho\downarrow0}\frac{E(B_\rho(a))}{\rho}=0
 \qquad\text{for a.e. interior axis point }a.
\]

For almost every $a\in\{q=P\}$, Lebesgue differentiation also gives
\[
 \frac{\HH^1(I\cap B_\rho(a))}{\rho}\longrightarrow0.
\]
Consequently,
\[
 \liminf_{\rho\downarrow0}\frac{\mu(B_\rho(a))}{\rho}
 =
 \liminf_{\rho\downarrow0}
 \left[
 \frac{E(B_\rho(a))}{\rho}
 +8\pi k\frac{\HH^1(I\cap B_\rho(a))}{\rho}
 \right]=0.
\]
Lemma~\ref{lem:small-energy}(2) therefore provides a neighborhood
of $a$ on which $u$ is regular and the normalized current vanishes.
Thus, locally on the axis,
\begin{equation}\label{eq:essential-closure}
 \HH^1\bigl(\{q=P\}\setminus\Reg(u)\bigr)=0.
\end{equation}

Now let $J$ be any relatively compact open interior axis interval.
If $\HH^1(J\setminus I)>0$, then \eqref{eq:essential-closure}
gives $J\cap\Reg(u)\ne\varnothing$.
Otherwise,
\[
 q=-P,\qquad m=P\qquad\text{a.e. on }J.
\]
Choose $a\in J$ with zero lower graph density and a sufficiently
small $\rho>0$ such that
\[
 B_\rho(a)\cap\axis\subset J,\qquad
 \frac{E(B_\rho(a))}{\rho}\le\epsilon_k.
\]
Lemma~\ref{lem:small-energy}(1) gives $a\in\Reg(u)$.
Hence every such $J$ contains a regular point.
\end{proof}

\section{Pole-level geometry and interior local finiteness}\label{sec:topology}

The planar level-set strategy is inherited from
\cite[proof of Theorem~9.1]{HLP}.

\subsection{Nodal graphs and finite axis contacts}\label{sec:nodal}

\begin{lemma}[Local nodal structure]\label{lem:nodal}
Let $T$ be locally minimizing, and suppose its graph map $u$ is
nonconstant on an axis-centered ball $B$ with $\overline B\Subset\Omega$.
On its open meridian half-disk $D$, set
\[
 Z_P=\{(r,z)\in D:u(r,0,z)=P\},\qquad P\in\{N,S\}.
\]
Each $Z_P$ is relatively closed and is a locally finite graph
of analytic curves. At every point, exactly $2d$ half-branches
meet for some integer $d\ge1$. In particular, there are no
isolated points or interior endpoints.
\end{lemma}

\begin{proof}
Relative closedness follows from off-axis continuity.
Fix $x_0\in Z_P$ and choose an analytic angle with
$\varphi(x_0)=j\pi$. Locally, $Z_P=\{f=0\}$ for
$f=\varphi-j\pi$, and
\[
 \Delta f+r^{-1}f_r-Vf=0,\qquad
 V=\frac{k^2}{r^2}\frac{\sin f\cos f}{f}\in C^\omega,
\]
where the quotient equals $1$ at zero.
Nonconstancy and analytic continuation give $f\not\equiv0$.
If $d\ge1$ is its vanishing order, then
\[
 f(x_0+y)=P_d(y)+O(|y|^{d+1}),\qquad \Delta P_d=0,
\]
so, for some $A\ne0$ and $t_0\in\R$,
\[
 P_d(\rho\cos t,\rho\sin t)
 =A\rho^d\cos(d(t-t_0)).
\]
The function
\[
 F(\rho,t)=\rho^{-d}f(x_0+\rho(\cos t,\sin t))
\]
extends analytically to $\rho=0$, where it has exactly $2d$
simple angular zeros. The implicit function theorem gives
$2d$ analytic half-branches; uniform convergence as $\rho\to0$
excludes additional zeros. Away from $x_0$ these branches are
regular and disjoint. A finite cover of any compact subset of
$D$ proves local finiteness.
\end{proof}

\begin{lemma}[Half-turn replacement]\label{lem:replacement}
A nonconstant minimizing graph has no pole-level cycle or embedded
pole-level arc with two interior axis limits, possibly equal, whose
Jordan interior is compactly contained in a minimizing ball.
\end{lemma}
\begin{proof}
Adjoin the intervening axis segment to the arc, or its common limit
when the limits coincide, and let $V$ be the Jordan interior.
For a cycle take its ordinary interior. Set
\[
 R_\pi=\operatorname{diag}(-1,-1,1),\qquad
 v=R_\pi u\text{ in }\widehat V,\quad v=u\text{ elsewhere},
 \qquad I_{T_v}=I_T.
\]
The two definitions agree on the off-axis interface. Segmentwise
Lipschitz estimates give local Lipschitz regularity there, and the
axis-extension argument gives $v\in H^1$. Moreover,
\[
 |\nabla v|=|\nabla u|\text{ a.e.},\qquad q_v=q_u,
 \qquad \partial T_v=0,\qquad \DD(T_v)=\DD(T).
\]
Lemma~\ref{lem:equality} gives $v=u$ off the axis. Hence
$u=R_\pi u$ on $\widehat V\setminus\axis$, so $u$ is a constant pole
on an open set. Analytic continuation contradicts nonconstancy.
\end{proof}


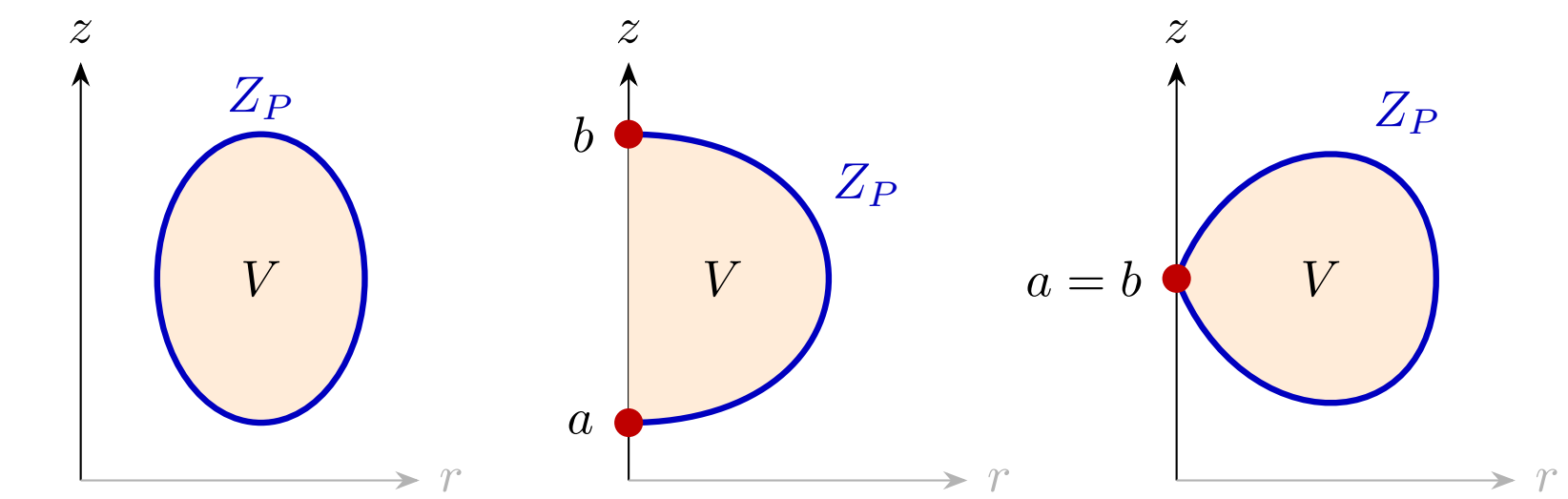
\begin{figure}[htbp]
\centering
\begin{tikzpicture}[
  font=\small,
  >=Stealth,
  level/.style={blue!75!black,very thick},
  region/.style={fill=orange!15},
  contact/.style={
    circle,fill=red!75!black,inner sep=2pt
  }
]

\foreach \shift in {0,3.8,7.6}{
  \begin{scope}[shift={(\shift,0)}]
    \draw[->] (0,-1.4)--(0,1.5) node[above] {$z$};
    \draw[->,gray!60] (0,-1.4)--(2.35,-1.4)
      node[right] {$r$};
  \end{scope}
}

\begin{scope}
  \draw[level,region]
    (1.25,0) ellipse [x radius=.72,y radius=1];
  \node at (1.25,0) {$V$};
  \node[blue!75!black] at (1.25,1.25) {$Z_P$};
  \node[align=center] at (1.1,-1.95)
    {(a) A pole-level cycle};
\end{scope}

\begin{scope}[shift={(3.8,0)}]
  \path[region]
    (0,-1)
    .. controls (1.85,-1) and (1.85,1) .. (0,1)
    -- cycle;
  \draw[level]
    (0,-1)
    .. controls (1.85,-1) and (1.85,1) .. (0,1);
  \node at (.65,0) {$V$};
  \node[blue!75!black] at (1.65,.65) {$Z_P$};
  \node[contact,label=left:{$a$}] at (0,-1) {};
  \node[contact,label=left:{$b$}] at (0,1) {};
  \node[align=center] at (1.1,-1.95)
    {(b) Two distinct\\axis limits};
\end{scope}

\begin{scope}[shift={(7.6,0)}]
  \draw[level,region]
    (0,0)
    .. controls (.45,1.15) and (1.8,1.15) .. (1.8,0)
    .. controls (1.8,-1.15) and (.45,-1.15) .. (0,0);
  \node at (1,0) {$V$};
  \node[blue!75!black] at (1.6,1.15) {$Z_P$};
  \node[contact,label=left:{$a=b$}] at (0,0) {};
  \node[align=center] at (1.1,-1.95)
    {(c) Two ends with\\the same axis limit};
\end{scope}

\node[align=center] at (4.85,-3.05) {
  $\displaystyle
  \widetilde u=
  \begin{cases}
    R_\pi u,&\text{in }\widehat V,\\
    u,&\text{outside }\widehat V,
  \end{cases}
  \qquad
  R_\pi=\operatorname{diag}(-1,-1,1),
  \qquad I_{\widetilde T}=I_T$
};

\node[align=center] at (4.85,-4.05) {
  $\displaystyle
  \DD(\widetilde T)=\DD(T)
  \ \Longrightarrow\ \widetilde u=u
  \ \Longrightarrow\
  u=R_\pi u\in\{N,S\}
  \quad\text{in }\widehat V\setminus\axis$
};

\end{tikzpicture}
\caption{The half-turn replacement excludes all three configurations.
Here $P\in\{N,S\}$, the blue curves lie in $Z_P$, and the shaded
Jordan region $V$ has its solid of revolution compactly contained
in the comparison ball.
Since $R_\pi$ fixes both poles, the replacement agrees across the
off-axis interface. It preserves the height, Dirichlet energy,
and vertical current. Equality and analytic continuation force
$u$ to be a constant pole on an open set, contradicting nonconstancy.
Red points indicate axis limits, not necessarily singularities.}
\label{fig:half-turn-replacement}
\end{figure}

Let $\mathcal G$ be the set of interior axis heights $h$ for which
$(0,0,h)$ is regular and, for some $\delta,\eta>0$,
\begin{equation}\label{eq:good-heights}
 u_1(r,0,z)\ne0
 \qquad(0<r<\delta,\ |z-h|<\eta).
\end{equation}
Thus no off-axis pole-level curve approaches an axis point with
height in $\mathcal G$.
This set is open and dense: on a regular analytic rectangle, write
\[
 u_1(r,0,z)=\sum_{j\ge0}A_j(z)r^j.
\]
For the first coefficient $A_{j_0}\not\equiv0$, its zeros are discrete,
and $A_{j_0}(h)\ne0$ gives
$|u_1(r,0,z)|\ge cr^{j_0}$ on a smaller rectangle.
Density of regular axis points then gives density of $\mathcal G$.

Fix an interior axis center $(0,t)$. Choose $s>0$ so that the
half-disk and its outer semicircle,
\[
 D_s=\{(r,z):r>0,\ r^2+(z-t)^2<s^2\},\qquad
 \Gamma_s=\{(r,z):r>0,\ r^2+(z-t)^2=s^2\},
\]
lie in the comparison neighborhood and $t\pm s\in\mathcal G$.
Sard's theorem allows $s$ to be chosen so that $\Gamma_s$ avoids
nodal vertices and meets every nodal edge transversely.
The endpoint neighborhoods contain no pole levels, so compactness
gives the finite outer-crossing set
\[
 \mathcal E=(Z_N\cup Z_S)\cap\Gamma_s,\qquad \#\mathcal E<\infty.
\]

Write
\[
 J_s=\{0\}\times(t-s,t+s),\qquad
 E=\overline{(Z_N\cup Z_S)\cap D_s}\cap J_s
\]
for the axis diameter and its pole-level contact set, respectively.
Consider a maximal graph path, following edges without backtracking.
It cannot repeat a vertex, since cycles are excluded by
Lemma~\ref{lem:replacement}.
If it does not exit through $\mathcal E$, local finiteness and
the good endpoint neighborhoods imply that its tail-limit set
\[
 K_\infty=\bigcap_{T>0}\overline{\gamma([T,\infty))}
\]
is nonempty, connected, and contained in $E$.
Since $E$ misses every height in the dense set $\mathcal G$,
it contains no nontrivial axis interval. Hence $K_\infty$ is a
singleton: every nonexiting path has a unique axis limit.

\begin{proposition}[Finite pole-level contacts and accessibility]
\label{prop:finite-poles}
In an enclosing half-disk chosen above, the axis contact set
$E=(\overline{Z_N}\cup\overline{Z_S})\cap\axis$ is finite, and
every contact is accessible by an embedded nodal arc.
\end{proposition}
\begin{proof}
Put $\mathcal Z_P=Z_P\cap D_s$ and
$\mathcal E_P=Z_P\cap\Gamma_s$. Fix $P$ and adjoin the finite
outer-crossing set $\mathcal E_P$ to $\mathcal Z_P$ as terminal vertices.
The graph is a forest. Every maximal path
has two ends, each an outer crossing or an axis limit. Two axis ends
are forbidden by Lemma~\ref{lem:replacement}, so every component $C$
meets $\mathcal E_P$ and
\[
 \#\pi_0(\mathcal Z_P)\le\#\mathcal E_P.
\]
Here $\pi_0(\mathcal Z_P)$ denotes its set of connected components. For $\mathcal E_P\cap C=\{b_1,\ldots,b_m\}$, let
$K_C=\bigcup_j[b_1,b_j]_C$ be the union of the unique graph paths.
These paths are compact and stay away from the axis, so local
finiteness makes $K_C$ a finite tree.
Every component of $C\setminus K_C$ attaches once to $K_C$ and
has a ray to the axis. Two remainders, or branching within one,
would yield a forbidden path with two axis ends. Thus
\[
 C=K_C\quad\text{or}\quad C=K_C\cup\gamma_C,
 \qquad\gamma_C\text{ is one unbranched ray with a unique axis limit}.
\]
There are finitely many components, each with at most one axis
contact, accessible along its ray. Therefore
\[
 \#(\overline{\mathcal Z_P}\cap\axis)\le\#\pi_0(\mathcal Z_P)\le\#\mathcal E_P.
\]
Summing over $P=N,S$ proves the assertion, including accumulation
points of whole components.
\end{proof}

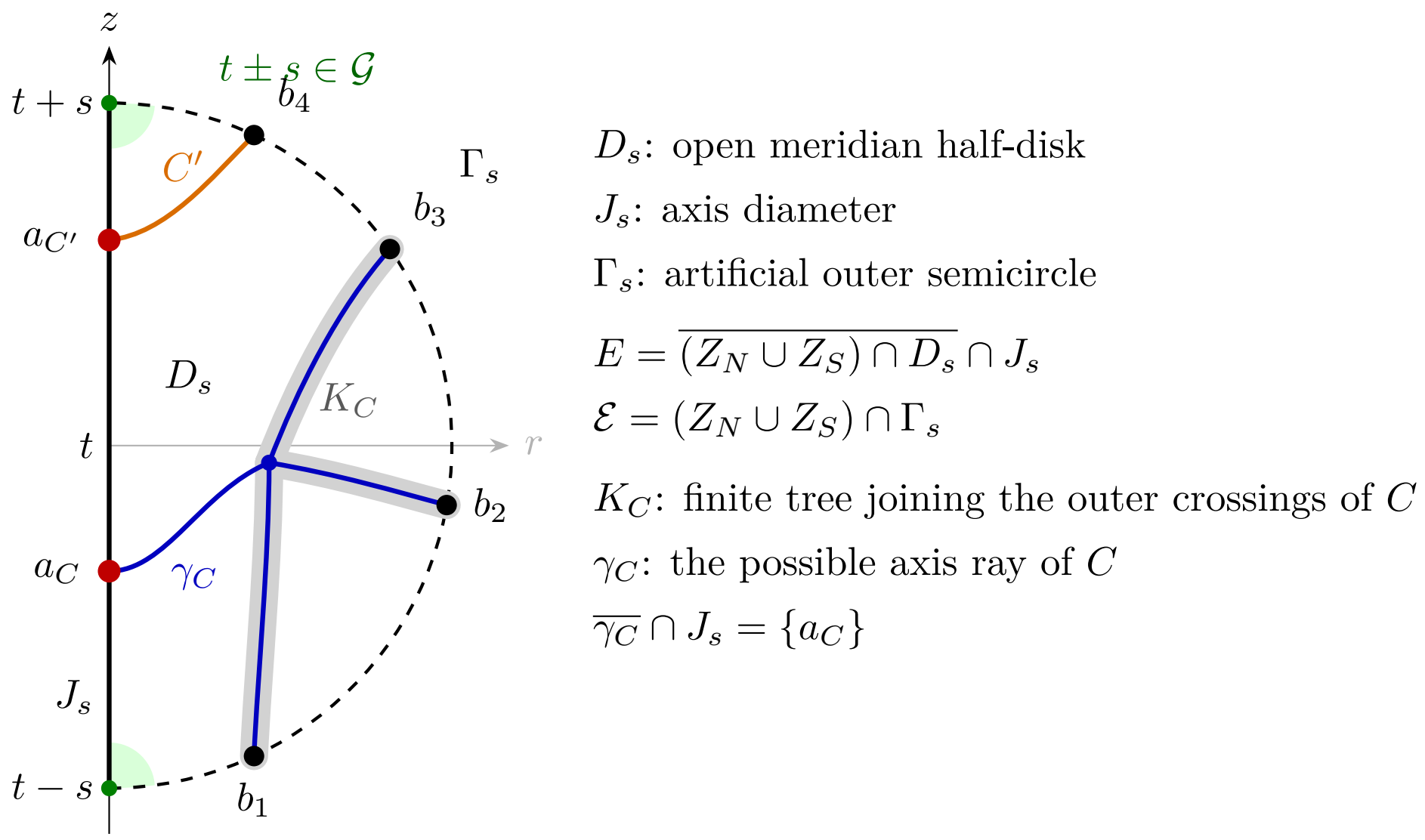
\begin{figure}[htbp]
\centering
\begin{tikzpicture}[
  >=Stealth,
  font=\small,
  level/.style={blue!75!black,very thick},
  core/.style={gray!35,line width=7pt,line cap=round},
  contact/.style={circle,fill=red!75!black,inner sep=2pt},
  crossing/.style={circle,fill=black,inner sep=1.8pt}
]

\begin{scope}
  \clip (0,-3)
    arc[start angle=-90,end angle=90,radius=3] -- cycle;
  \fill[green!15] (0,3) circle (.4);
  \fill[green!15] (0,-3) circle (.4);
\end{scope}

\draw[->,gray!60] (0,0) -- (3.5,0) node[right] {$r$};
\draw[->] (0,-3.4) -- (0,3.5) node[above] {$z$};
\draw[very thick] (0,-3) -- (0,3);
\draw[thick,dashed]
  (0,-3) arc[start angle=-90,end angle=90,radius=3];

\node[left] at (0,3) {$t+s$};
\node[left] at (0,-3) {$t-s$};
\node[left] at (0,0) {$t$};
\node[left] at (0,-2.2) {$J_s$};
\node at (.7,.6) {$D_s$};
\node at (3.25,2.45) {$\Gamma_s$};

\fill[green!50!black] (0,3) circle (2pt);
\fill[green!50!black] (0,-3) circle (2pt);
\node[green!40!black] at (1.65,3.3)
  {$t\pm s\in\mathcal G$};

\coordinate (v) at (1.4,-.15);
\coordinate (a) at (0,-1.1);
\coordinate (b1) at (-65:3);
\coordinate (b2) at (-10:3);
\coordinate (b3) at (35:3);

\draw[core]
  (v) .. controls (1.4,-1) and (1.3,-2) .. (b1);
\draw[core]
  (v) .. controls (2,-.25) and (2.5,-.4) .. (b2);
\draw[core]
  (v) .. controls (1.65,.5) and (2,1.2) .. (b3);

\draw[level]
  (v) .. controls (1.4,-1) and (1.3,-2) .. (b1);
\draw[level]
  (v) .. controls (2,-.25) and (2.5,-.4) .. (b2);
\draw[level]
  (v) .. controls (1.65,.5) and (2,1.2) .. (b3);

\draw[level]
  (a) .. controls (.5,-1.1) and (.8,-.4) .. (v);
\node[blue!75!black] at (.75,-1.15) {$\gamma_C$};
\node[gray!70!black] at (2.1,.4) {$K_C$};
\fill[blue!75!black] (v) circle (2pt);

\coordinate (aa) at (0,1.8);
\coordinate (b4) at (65:3);
\draw[orange!85!black,very thick]
  (aa) .. controls (.5,1.85) and (.9,2.35) .. (b4);
\node[orange!85!black] at (.65,2.45) {$C'$};

\node[contact,label=left:{$a_C$}] at (a) {};
\node[contact,label=left:{$a_{C'}$}] at (aa) {};

\node[crossing,label=below:{$b_1$}] at (b1) {};
\node[crossing,label=right:{$b_2$}] at (b2) {};
\node[crossing,label=above right:{$b_3$}] at (b3) {};
\node[crossing,label=above right:{$b_4$}] at (b4) {};

\node[anchor=north west,align=left] at (4.1,2.9) {
  $D_s$: open meridian half-disk\\[4pt]
  $J_s$: axis diameter\\[4pt]
  $\Gamma_s$: artificial outer semicircle\\[8pt]
  $\displaystyle
    E=\overline{(Z_N\cup Z_S)\cap D_s}\cap J_s$\\[4pt]
  $\displaystyle
    \mathcal E=(Z_N\cup Z_S)\cap\Gamma_s$\\[8pt]
  $K_C$: finite tree joining the outer crossings of $C$\\[4pt]
  $\gamma_C$: the possible axis ray of $C$\\[4pt]
  $\displaystyle
    \overline{\gamma_C}\cap J_s=\{a_C\}$
};

\end{tikzpicture}
\caption{Finite axis contacts from finitely many outer crossings.
After adjoining the outer-crossing vertices, each pole-level
component consists of a finite core and at most one axis ray.
Red points are axis contacts; black points are outer crossings.
Axis contacts need not be singular points of the map.}
\label{fig:finite-axis-contacts}
\end{figure}
\FloatBarrier

\subsection{At most two singularities on a one-sided axis interval}
\label{sec:comparisons}

\begin{proposition}
\label{prop:two-points}
Let $T$ be locally minimizing, and let
$C=(0,R)\times I$, where $I$ is a bounded open interval and
$\widehat C\Subset\Omega$. Suppose an angle function satisfies
\begin{equation}\label{eq:one-sided}
 0<\varphi(r,z)<\pi\qquad\text{on }C.
\end{equation}
Then
\[
 \#\{z\in I:(0,0,z)\in\Sing(u)\}\le2.
\]
\end{proposition}

\begin{proof}
Reflecting the target height and replacing $I_T$ by $-I_T$ if
necessary, assume $c<0$. The lower-current pole is then $S$.

\smallskip\noindent
\textit{Step 1: pole estimates and removal.}
At a regular axis point, equivariance permits only transverse
Taylor monomials $\zeta^p\overline\zeta^{\,\ell}$ with
$p-\ell=k$, where $\zeta=x_1+ix_2$. Hence, uniformly on a smaller
axis neighborhood,
\begin{equation}\label{eq:pole-order}
 |u_1+iu_2|\le Cr^k,\qquad
 \begin{cases}
 \varphi\le Cr^k,&\text{near a regular north point},\\
 \pi-\varphi\le Cr^k,&\text{near a regular south point}.
 \end{cases}
\end{equation}

Conversely, suppose $|\varphi|\le Cr^k$ near an axis interval.
For the physical scalar $\Phi(x)=\varphi(r,z)$,
\[
 \Delta\Phi=k^2r^{-2}\sin\Phi\cos\Phi.
\]
Interior Poisson estimates on balls of radius comparable to $r$
give
\[
 |\nabla\varphi|
 \le C\bigl(r^{-1}\|\Phi\|_\infty
             +r\|\Delta\Phi\|_\infty\bigr)
 \le Cr^{k-1},
 \qquad
 \frac{k|\sin\varphi|}{r}\le Cr^{k-1}.
\]
Thus $\nabla u$ is locally bounded. Equation~\eqref{eq:weak-harmonic},
interior $W^{2,p}$ estimates, and bootstrapping give smoothness
across the axis, with north value. The same argument applies at
the south pole when $|\pi-\varphi|\le Cr^k$.

\smallskip\noindent
\textit{Step 2: square completions.}
For $C_0=(0,\rho)\times(a,b)$ with
$\overline{\widehat C_0}\subset\widehat C$, define
\[
 Q_\pm(\varphi)
 =2\pi\int_{C_0}r
 \left[
 \left(\varphi_r\mp\frac{k\sin\varphi}{r}\right)^2
 +\varphi_z^2
 \right]dr\,dz.
\]
Using $m=(c-q)/2$, $c<0$, and the axis trace,
\begin{align}
 \DD(T;\widehat C_0)
 &=Q_-(\varphi)
   +4\pi k\int_a^b
       \bigl(\cos\varphi(\rho,z)-c\bigr)\,dz,
 \label{eq:qminus}\\
 E(u;\widehat C_0)
 &=Q_+(\varphi)
   +4\pi k\int_a^b
       \bigl(q(z)-\cos\varphi(\rho,z)\bigr)\,dz.
 \label{eq:qplus}
\end{align}
The first identity permits changes of the axis trace and current;
the second will be used with both traces fixed.

For a fixed $\epsilon>0$, set
\begin{equation}\label{eq:radial-calibration}
 \beta_\epsilon(r)=2\arctan\frac{\epsilon}{r^k},
 \qquad
 \xi_\epsilon(r)=2\arctan\frac{r^k}{\epsilon}.
\end{equation}
They satisfy
\[
 (\beta_\epsilon)_r+\frac{k\sin\beta_\epsilon}{r}=0,
 \qquad
 (\xi_\epsilon)_r-\frac{k\sin\xi_\epsilon}{r}=0,
 \qquad
 Q_-(\beta_\epsilon)=Q_+(\xi_\epsilon)=0.
\]
Each barrier has energy at most $8\pi k(b-a)$.
The max/min clamps below therefore have finite energy and extend
across the axis in $H^1$.

\smallskip\noindent
\textit{Step 3: two clamping comparisons.}
First suppose the cap axis points at heights $a,b$ are regular
south points. By \eqref{eq:pole-order} and compactness away from
the axis, sufficiently small $\epsilon$ gives
\[
 r^k\tan(\varphi/2)>\epsilon,
 \qquad\text{hence }\varphi>\beta_\epsilon,
\]
on full boundary collars of a suitably narrow cylinder.
Set
\[
 \widetilde\varphi=\max\{\varphi,\beta_\epsilon\},
 \qquad \widetilde q=-1,\qquad
 \widetilde m=\frac{c+1}{2}
\]
inside the cylinder, retaining the original map and current outside.
They agree on the collars, and
$\widetilde q+2\widetilde m=c$ proves closure.
Equation~\eqref{eq:qminus} gives
\[
 \DD(T;\widehat C_0)-\DD(\widetilde T;\widehat C_0)
 =
 2\pi\int_{\{\varphi<\beta_\epsilon\}}r
 \left[
 \left(\varphi_r+\frac{k\sin\varphi}{r}\right)^2
 +\varphi_z^2
 \right]dr\,dz\ge0.
\]
Lemma~\ref{lem:equality} and collar agreement give
$u=\widetilde u$, hence
\[
 \varphi\ge\beta_\epsilon,\qquad
 \pi-\varphi\le2\arctan(r^k/\epsilon)\le2r^k/\epsilon.
\]
Step~1 makes the entire intervening axis interval regular south.

Next suppose $q=1$ almost everywhere on $(a,b)$ and both cap
axis points are regular north. Equation~\eqref{eq:pole-order}
and compactness give, for small $\epsilon$,
\[
 r^{-k}\tan(\varphi/2)<\epsilon^{-1},
 \qquad\text{hence }\varphi<\xi_\epsilon
\]
on full boundary collars. Set
\[
 \widetilde\varphi=\min\{\varphi,\xi_\epsilon\},
 \qquad \widetilde q=q=1,\qquad \widetilde m=m.
\]
The current cost is unchanged, and \eqref{eq:qplus} gives
\[
 \DD(T;\widehat C_0)-\DD(\widetilde T;\widehat C_0)
 =
 2\pi\int_{\{\varphi>\xi_\epsilon\}}r
 \left[
 \left(\varphi_r-\frac{k\sin\varphi}{r}\right)^2
 +\varphi_z^2
 \right]dr\,dz\ge0.
\]
Again Lemma~\ref{lem:equality} yields
\[
 \varphi\le\xi_\epsilon\le2r^k/\epsilon.
\]
Step~1 makes the entire interval regular north.
The assumption $q=1$ is essential here to retain the original current.

\smallskip\noindent
\textit{Step 4: count the possible singular points.}
Let
\[
 V_S=\{z\in I:(0,0,z)\in\Reg(u),\ u(0,0,z)=S\}.
\]
This set is open, and the south comparison gives
\[
 a,b\in V_S,\ a<b\quad\Longrightarrow\quad(a,b)\subset V_S.
\]
Thus $V_S$ is an interval, possibly empty or all of $I$.
Since the normalized lower-current pole is $S$,
\eqref{eq:essential-closure} gives
\[
 \HH^1(\{q=-1\}\setminus V_S)=0.
\]
Consequently $q=1$ almost everywhere on each component
$J$ of $I\setminus\overline{V_S}$.
By Lemma~\ref{lem:dense}, every $z\in J$ can be bracketed by
regular north heights $a<z<b$ in $J$.
The north comparison makes $z$ regular. Therefore
\[
 \{z\in I:(0,0,z)\in\Sing(u)\}\subset\partial_I V_S,
 \qquad \#\partial_I V_S\le2.
\]
\end{proof}

\subsection{Interior local finiteness}\label{sec:assembly}

We combine the finite contact count from
Proposition~\ref{prop:finite-poles} with the two-point bound from
Proposition~\ref{prop:two-points}. 

\begin{lemma}[Interior local finiteness]\label{lem:interior}
Let $k\ge2$ and let $T\in\operatorname{cart}_{\mathrm{AS}}^{(k)}(\Omega,\Sph)$ be locally minimizing.
Then
\[
 \Sing(u)\subset\axis,\qquad
 \#(\Sing(u)\cap K)<\infty
 \quad\text{for every compact }K\Subset\Omega.
\]
\end{lemma}

\begin{proof}
Fix an interior axis point. The construction in
Subsection~\ref{sec:nodal} gives
an enclosing half-disk $D$ with axis diameter $I$.
The constant-map case is immediate. Otherwise let
\[
 E=(\overline{Z_N}\cup\overline{Z_S})\cap I,\qquad e=\#E<\infty
\]
by Proposition~\ref{prop:finite-poles}. For a component $J$ of
$I\setminus E$ and $K\Subset J$, compactness gives $R_K>0$ and
a slightly larger interval $K'\Supset K$ with
\[
 ((0,R_K)\times K')\cap(Z_N\cup Z_S)=\varnothing.
\]
Indeed, failure would give $r_n\downarrow0$, $z_n\to z\in K$ and
$(r_n,z_n)\in Z_N\cup Z_S$, whence $z\in E$, a contradiction.
A real lift on this simply connected rectangle satisfies
\[
 \varphi((0,R_K)\times K')\subset(j\pi,(j+1)\pi)
\]
for one integer $j$. If $j$ is even subtract $j\pi$. If $j$ is
odd apply the antipodal target isometry and simultaneously send
$I_T$ to $-I_T$, then subtract the appropriate multiple of $2\pi$.
These simultaneous transformations of the full current put the
angle in $(0,\pi)$ without changing the energy or minimality.
Proposition~\ref{prop:two-points} gives
\[
 \#(\Sing(u)\cap K)\le2
 \quad\Longrightarrow\quad \#(\Sing(u)\cap J)\le2,
\]
since any three distinct points of $J$ lie in one compact subinterval.
Consequently
\[
 \#(\Sing(u)\cap\widehat D)\le e+2\#\pi_0(I\setminus E)
 \le e+2(e+1)=3e+2.
\]
There are no off-axis singularities by
Subsection~\ref{sec:scalar-equation}.
For $K_0\Subset\Omega$ choose a finite cover $\{\widehat D_i\}_{i=1}^m$
of $K_0\cap\axis$ by such neighborhoods. Then
\[
 \#(\Sing(u)\cap K_0)\le\sum_{i=1}^m(3e_i+2)<\infty.
\]
\end{proof}

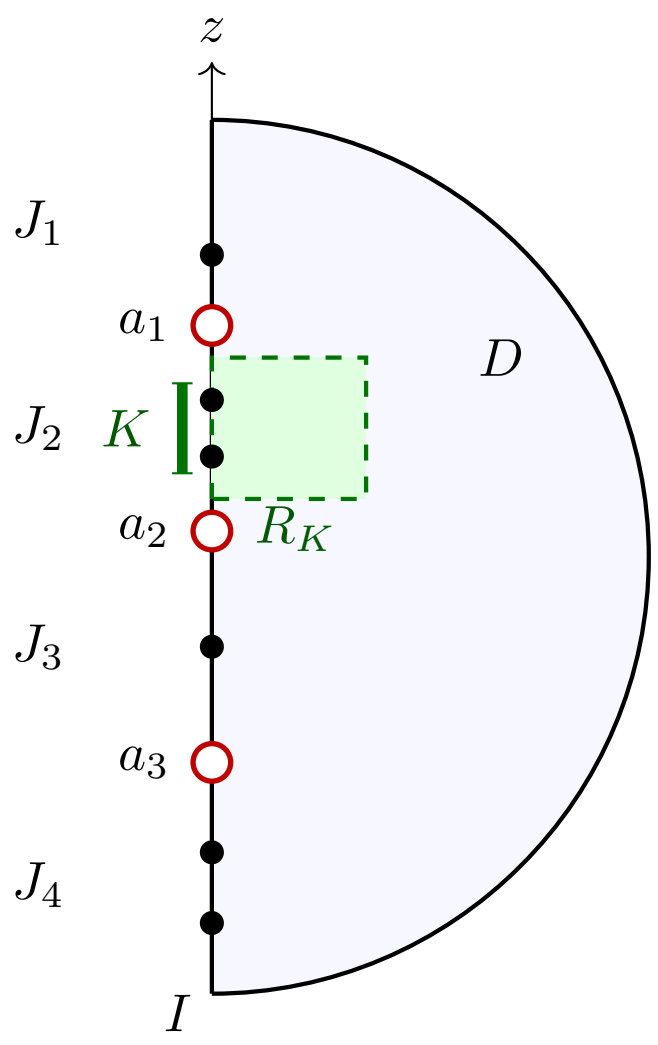
\begin{figure}[htbp]
\centering

\begin{minipage}[c]{0.43\textwidth}
\centering
\begin{tikzpicture}[
  x=.85cm,y=.85cm,
  font=\small,
  contact/.style={
    circle,draw=red!75!black,fill=white,
    line width=1pt,inner sep=2.5pt
  },
  singular/.style={
    circle,fill=black,inner sep=1.6pt
  }
]

\fill[blue!3]
  (0,-3.4)
  arc[start angle=-90,end angle=90,radius=3.4]
  -- cycle;
\draw[thick]
  (0,-3.4)
  arc[start angle=-90,end angle=90,radius=3.4];
\draw[thick] (0,-3.4)--(0,3.4);
\draw[->] (0,3.4)--(0,3.85) node[above] {$z$};

\node at (2.25,1.55) {$D$};
\node[left] at (0,-3.55) {$I$};

\fill[green!12] (0,.45) rectangle (1.2,1.55);
\draw[green!45!black,dashed,thick]
  (0,.45) rectangle (1.2,1.55);
\node[green!35!black] at (.65,.22) {$R_K$};

\draw[green!45!black,line width=2pt]
  (-.23,.65)--(-.23,1.35);
\draw[green!45!black]
  (-.31,.65)--(-.15,.65)
  (-.31,1.35)--(-.15,1.35);
\node[left,green!35!black] at (-.32,1) {$K$};

\foreach \h in {2.35,.78,1.22,-.70,-2.30,-2.85}
  \node[singular] at (0,\h) {};

\node[contact,label=left:{$a_1$}] at (0,1.8) {};
\node[contact,label=left:{$a_2$}] at (0,.2) {};
\node[contact,label=left:{$a_3$}] at (0,-1.6) {};

\node at (-1.35,2.60) {$J_1$};
\node at (-1.35,1.00) {$J_2$};
\node at (-1.35,-.70) {$J_3$};
\node at (-1.35,-2.55) {$J_4$};

\end{tikzpicture}
\end{minipage}
\hfill
\begin{minipage}[c]{0.53\textwidth}
\small
\setlength{\parindent}{0pt}

\textbf{1. Finite contact set.}
\[
 E=(\overline{Z_N}\cup\overline{Z_S})\cap I,
 \qquad e=\#E<\infty.
\]
Thus
\[
 \#\pi_0(I\setminus E)\le e+1.
\]

\medskip
\textbf{2. One-sided comparison.}

For every compact interval $K\Subset J_j$, choose
$K\Subset K'\Subset J_j$ and $R_K>0$ such that
\[
 \bigl((0,R_K)\times K'\bigr)
 \cap(Z_N\cup Z_S)=\varnothing.
\]
After normalization, $0<\varphi<\pi$. Hence
\[
 \#(\Sing(u)\cap K)\le2.
\]
Exhausting $J_j$ gives
\[
 \#(\Sing(u)\cap J_j)\le2.
\]

\medskip
\textbf{3. Counting.}
\[
 \begin{aligned}
 \#(\Sing(u)\cap\widehat D)
 &\le e+2\#\pi_0(I\setminus E)\\
 &\le 3e+2<\infty.
 \end{aligned}
\]

\end{minipage}

\caption{Interior local finiteness.
Red circles mark the finite axis contact set $E$; they need not be
singular. Black dots illustrate possible singularities, with at most
two in each component $J_j$ of $I\setminus E$.
The shaded rectangle is a pole-free neighborhood above
$K\Subset J_2$. A finite cover by such half-disks proves finiteness
on every compact subset of the domain.}
\label{fig:interior-local-finiteness}
\end{figure}
\FloatBarrier

\section{The map--current correspondence}\label{sec:map}

Following Hardt--Lin--Poon \cite[Sections~4--5]{HLP}, we identify
the canonical axis current for degree-zero boundary data and show
that its energy equals $E_1$. We then prove existence and transfer
local map minimality to the structural current problem.

\subsection{Boundary data and the canonical connection}

Let $g\in H^1(\partial B_R,\Sph)$ be meridional $k$-axial.
Its angle profile satisfies
\begin{equation}\label{eq:boundary-profile-energy}
 \int_{\partial B_R}|\nabla_{\rm tan}g|^2
 =2\pi\int_0^\pi
 \left(\sin\alpha\,|\varphi_g'|^2
       +\frac{k^2\sin^2\varphi_g}{\sin\alpha}\right)d\alpha<\infty.
\end{equation}
Cauchy--Schwarz gives $(\cos\varphi_g)'\in L^1(0,\pi)$.
The endpoint limits therefore exist and, by finite energy, belong
to $\{-1,1\}$. At either endpoint with limiting height $P$,
\[
 |g-(0,0,P)|^2=2(1-P\cos\varphi_g)
\]
gives uniform convergence in longitude.
Together with one-dimensional $H^1$ continuity away from the poles,
this gives a continuous representative of $g$, which we use throughout.

For degree-zero data, the following formula identifies the canonical
axis current and its mass, extending \cite[Lemma~4.1]{HLP} to weight $k$.

\begin{lemma}\label{lem:L-length}
Let $u\in\mathcal A_k(g)$, with meridional $H^1$ degree-zero datum $g$.
Its two boundary poles have a common height $P\in\{-1,1\}$, and
\begin{equation}\label{eq:L-axis}
 L(u;B_R)=\frac{k}{2}\int_{-R}^R|P-q(z)|\,dz.
\end{equation}
The canonical current
\begin{equation}\label{eq:canonical-current}
 I_u=\frac{k}{2}(P-q)\curr{\axis\cap B_R}
\end{equation}
has mass $L(u;B_R)$ and multiplicities $0$ and $Pk$.
\end{lemma}
\begin{proof}
The Sobolev degree formula gives
\[
 0=\deg g
 =-\frac{k}{2}\int_0^\pi(\cos\varphi_g)'\,d\alpha
 =\frac{k}{2}(P_+-P_-),
\]
so $P_+=P_-=P\in\{-1,1\}$.

Let $\xi$ be a smooth test function with $|\nabla\xi|\le1$, let
$\bar\xi$ be its angular average, and put $\psi(z)=\xi(0,0,z)$.
Orient the physical semicircle $\gamma$ from bottom to top.
Since $\cos\varphi\in W^{1,1}(D_R)$ has traces $q$ on the
axis and $\cos\varphi_g$ on $\gamma$, Gauss--Green yields
\[
 \begin{aligned}
 &\int_{B_R}D(u)\cdot\nabla\xi
 -\int_{\partial B_R}(\operatorname{Jac}g)\xi\\
 &\qquad=2\pi k\left[
 \int_{-R}^R q\psi'\,dz
 -\int_\gamma d(\cos\varphi_g\,\bar\xi)\right]\\
 &\qquad=2\pi k\left[
 \int_{-R}^R q\psi'\,dz-P\bigl(\psi(R)-\psi(-R)\bigr)\right].
 \end{aligned}
\]
Consequently,
\begin{equation}\label{eq:L-dual-axis}
 \frac1{4\pi}\left[
 \int_{B_R}D(u)\cdot\nabla\xi
 -\int_{\partial B_R}(\operatorname{Jac}g)\xi\right]
 =\frac{k}{2}\int_{-R}^R(q-P)\psi'\,dz.
\end{equation}
Since $|\psi'|\le1$,
\[
 L(u;B_R)\le\frac{k}{2}\int_{-R}^R|q-P|\,dz.
\]
Equality is attained by $\xi(x)=-Px_3$, because
\[
 (q-P)(-P)=1-Pq=|q-P|.
\]
Finally,
\[
 \mathbf M(I_u)
 =\frac{k}{2}\int_{-R}^R|P-q|\,dz
 =L(u;B_R),
\]
and the coefficient of $I_u$ belongs to $\{0,Pk\}$.
\end{proof}

\subsection{Existence and local minimality}

Only the canonical compatibility constant $c=P$ is needed for the
fixed-boundary problem. We first prove existence without assuming
convergence of the axis traces, and then transfer minimality to the
structural current class.

\begin{lemma}[Existence]\label{lem:existence}
For every meridional $k$-axial datum
$g\in H^1(\partial B_R,\Sph)$ of degree zero,
$E_1$ attains its minimum in $\mathcal A_k(g)$.
\end{lemma}
\begin{proof}
Let $D_R$ be the meridian half-disk, $\Gamma_R$ its physical
semicircle, and $P$ the common boundary pole height.
For the circle profile $w=(w_1,w_2)=(\sin\varphi,\cos\varphi)$, put
$j=w_2\nabla w_1-w_1\nabla w_2=\nabla\varphi$.
Gauss--Green and Lemma~\ref{lem:L-length} give
\begin{equation}\label{eq:canonical-square}
 \begin{split}
 Q_P(u)&=2\pi\int_{D_R}r
 \left[\left(j_r-\frac{Pk w_1}{r}\right)^2+j_z^2\right],\\
 E_1(u)&=Q_P(u)+4\pi k
 \left(2R-P\int_{\Gamma_R}(\cos\varphi_g)\nu_r\,ds\right).
 \end{split}
\end{equation}
The last term is fixed by $g$.

The radial extension $u_0(x)=g(Rx/|x|)$ satisfies
\[
 E(u_0)=R\int_{\partial B_R}|\nabla_{\rm tan}g|^2<\infty,
 \qquad L(u_0)\le2kR.
\]
A minimizing sequence therefore has a subsequence satisfying
\[
 u_n\rightharpoonup u\text{ in }H^1,\qquad
 u_n\to u\text{ in }L^2\text{ and a.e.},\qquad
 u\in\mathcal A_k(g).
\]
On $D_{R,\delta}=D_R\cap\{r>\delta\}$, the circle profiles satisfy
\[
 w_n\to w\text{ a.e.},\qquad |w_n|=1,\qquad
 j_n\rightharpoonup j\text{ in }L^2,
 \qquad w_{n,1}/r\to w_1/r\text{ in }L^2.
\]
Here weak convergence of $j_n$ follows by testing its expression
in terms of $w_n$ and $\nabla w_n$; boundedness and a.e. convergence
make the profile factors converge strongly against every $L^2$ test.
Weak lower semicontinuity, followed by $\delta\downarrow0$, yields
\[
 Q_P(u)\le\liminf_n Q_P(u_n),\qquad
 E_1(u)\le\liminf_n E_1(u_n)=\inf_{\mathcal A_k(g)}E_1.
\]
\end{proof}

The next lemma is the localization step of
\cite[Section~5.2]{HLP}, with the weight-$k$ connection retained.

\begin{lemma}[The canonical lift]\label{lem:canonical-lift}
Let $u\in\mathcal A_k(g)$, where $g$ is meridional $k$-axial
of degree zero, and define
\[
 T_u=\graph{u}+I_u\times\curr{\Sph},\qquad
 I_u=\frac{k}{2}(P-q)\curr{\axis\cap B_R}.
\]
Then $T_u\in\operatorname{cart}_{\mathrm{AS}}^{(k)}(B_R,\Sph)$ and
\begin{equation}\label{eq:canonical-energy}
 \DD(T_u;B_R)=E_1(u;B_R).
\end{equation}
If $u$ locally minimizes $E_1$ under compactly supported meridional
changes, then $T_u$ is locally minimizing.
\end{lemma}
\begin{proof}
Since $q+2(P-q)/2=P$, Lemma~\ref{lem:full-boundary} gives
$\partial T_u=0$ in $B_R\times\Sph$.
Lemma~\ref{lem:L-length} gives \eqref{eq:canonical-energy}.

Let $\widetilde T=\graph{v}+km_v\curr{\axis}\times\curr{\Sph}$
be an admissible compact perturbation of $T_u$.
Then $v=u$ outside a compact subset of $B_R$, so $v\in\mathcal A_k(g)$.
Closure gives $q_v+2m_v=c_v$, while exterior agreement on the
axis fixes $c_v=P$. Hence
\[
 \widetilde T=T_v,\qquad
 \DD(T_u;B_R)=E_1(u)\le E_1(v)=\DD(\widetilde T;B_R).
\]
Cancelling the unchanged exterior energy proves local minimality.
\end{proof}

\section{Global finiteness for continuous degree-zero data}
\label{sec:boundary}

Let $u$ minimize $E_1$ in $\mathcal A_k(g)$ under the assumptions
of Theorem~\ref{thm:global}. Dilation and a target height reflection
reduce the proof to $R=1$ and
\begin{equation}\label{eq:boundary-canonical}
 g(p_+)=g(p_-)=N,\qquad
 E_1(u)=E(u)+8\pi k\int_{-1}^1\mathbf1_{\{q=-1\}}\,dz,
 \qquad p_\pm=(0,0,\pm1).
\end{equation}
Lemmas~\ref{lem:canonical-lift} and~\ref{lem:interior} give
$\Sing(u)\subset\axis$ and interior local finiteness.
It remains to control the two axis tails near $p_\pm$.
We treat $p_+$ and assume $u$ is nonconstant.
Write $D_1=\{(r,z):r>0,\ r^2+z^2<1\}$, put $p=(0,1)$,
and use the pole-level sets
\[
 Z_P=\{(r,z)\in D_1:u(r,0,z)=P\},\qquad P\in\{N,S\},
\]
as in Subsection~\ref{sec:nodal}. In meridian expressions, $u$ also
denotes its restriction $(r,z)\mapsto u(r,0,z)$.

\subsection{Boundary continuity and comparison}

\begin{lemma}[Off-axis boundary continuity]\label{lem:off-axis-boundary}
The representative of $u$ is continuous on
$\overline B_1\setminus\axis$ and equals $g$ on
$\partial B_1\setminus\axis$.
\end{lemma}
\begin{proof}
Fix an off-axis boundary arc, where $r\ge r_*>0$.
Extend its circle profile across the arc and apply the lifting
construction of Subsection~\ref{sec:scalar-equation}.
This gives an $H^1$ angle $\varphi$ with continuous $H^1$ boundary
angle $\varphi_g$, after an additive adjustment by $2\pi\Z$.
For $\Psi=\sqrt r\,\varphi$, equation~\eqref{eq:pde} gives
\[
 \Delta\Psi=k^2r^{-3/2}\sin\varphi\cos\varphi
              -\frac{\Psi}{4r^2}=:F,
 \qquad F\in L^p\quad(1<p<\infty)
\]
locally up to the arc. Subtract a Newton potential
$W\in W^{2,p}$, $p>2$, with $\Delta W=F$.
The harmonic function $H=\Psi-W$ has continuous boundary trace.
After conformally straightening the arc, write
\[
 H=\mathsf P b+h,\qquad
 b\in H^{1/2}(\R)\cap C_c^0(\R),\qquad
 \operatorname{Tr}h=0
\]
on a smaller arc, where $\mathsf P b$ is the Poisson extension
of a cutoff of the trace. Then
$\mathsf P b(\cdot,t)\to b$ locally uniformly, while the odd
reflection of $h$ is harmonic. Thus $H$, $\varphi$, and $u$
are continuous up to the arc with the prescribed traces.
\end{proof}

The replacements of Lemmas~\ref{lem:equality} and~\ref{lem:replacement}
remain valid for the Jordan regions $V\subset D_1$ used below,
whose closures meet the physical semicircle only at $p$.
Indeed, let $A:\Sph\to\Sph$ be a Lipschitz target map preserving
the meridional class and fixing $u$ on the off-axis interface.
Set $v=A\circ u$ in $\widehat V$ and $v=u$ elsewhere.
The same off-axis pasting and axis-extension arguments give
$v\in H^1(B_1,\Sph)$ whenever the Dirichlet energy does not increase.
Since $v=u$ near every physical boundary point except possibly
$p_+$, its Sobolev trace is still $g$. Consequently,
\begin{equation}\label{eq:boundary-equality}
 E(v)\le E(u),\quad L(v)\le L(u)
 \quad\Longrightarrow\quad
 E_1(u)\le E_1(v)\le E_1(u),\qquad v=u\text{ off the axis}.
\end{equation}
The last equality follows, as in Lemma~\ref{lem:equality}, from
analytic continuation and agreement on an exterior open set.
All interfaces below are locally rectifiable away from the axis,
so the pasting argument applies.

In particular, take $A=R_\pi=\operatorname{diag}(-1,-1,1)$.
It fixes both poles and preserves the height and both energy terms:
\begin{equation}\label{eq:boundary-half-turn}
 q_v=q,\qquad E_1(v)=E_1(u),\qquad
 u=v=R_\pi u\quad\text{in }\widehat V\setminus\axis.
\end{equation}
As in Lemma~\ref{lem:replacement}, this would make $u$ a constant
pole on an open set. Hence pole-level cycles and arcs with two
closed-axis ends are excluded whenever their Jordan interiors
meet the physical boundary only at $p$; one end may now be $p$
itself, and the two ends may coincide.

\subsection{Finiteness near the boundary poles}

\begin{proposition}[No singular accumulation at a boundary pole]
\label{prop:boundary-nonaccumulation}
There exists $t_*<1$ such that
\[
 \#\{z\in(t_*,1):(0,0,z)\in\Sing(u)\}<\infty.
\]
\end{proposition}
\begin{proof}
\textit{Step 1: a cap collar and a finite interior section.}
Fix $0<\alpha_0<\alpha\le\pi/4$.
Lemma~\ref{lem:off-axis-boundary} and $g(p_+)=N$ give a relative
open collar of the physical arc near $p$ whose interior image
lies in $B_{\alpha_0}(N)$.
With $\rho(z)=\sqrt{1-z^2}$, choose smooth functions $a,b$ on
$(z_*,1)$, for some $z_*<1$, such that
\begin{equation}\label{eq:collar-walls}
 0<a(z)<b(z)<\rho(z),\qquad
 u(r,0,z)\in B_{\alpha_0}(N)\quad\text{if }a(z)\le r\le b(z).
\end{equation}
The walls can be chosen by a smooth minorant of the collar thickness.
Since $b<\rho$, both tend to zero as $z\uparrow1$.
Choose $z_0\in(z_*,1)\cap\mathcal G$ using
\eqref{eq:good-heights}, and perturb $a(z_0)$ within the collar
so that $u_1(a(z_0),0,z_0)\ne0$.
Set
\[
 \begin{split}
 J&=(z_0,1),\qquad
 \Gamma_0=\{(r,z_0):0<r<a(z_0)\},\\
 D_a&=\{(r,z):z\in J,\ 0<r<a(z)\},\qquad
 D_b=\{(r,z):z\in J,\ 0<r<b(z)\}.
 \end{split}
\]
The function $r\mapsto u_1(r,0,z_0)$ is analytic and is nonzero
near both endpoints of $\Gamma_0$, by the good-height condition
and the choice of $a(z_0)$. Consequently,
\begin{equation}\label{eq:bottom-ports}
 \mathcal E_0=\Gamma_0\cap(Z_N\cup Z_S),\qquad
 M:=\sum_{x\in\mathcal E_0}2d_x<\infty,
\end{equation}
where $2d_x$ is the number of incident nodal half-branches at $x$
from Lemma~\ref{lem:nodal}. Thus tangencies are counted as well.

Let $\mathcal U$ be the connected component of
$D_b\cap u^{-1}(B_{\alpha_0}(N))$ containing the connected strip
$\{(r,z):z\in J,\ a(z)\le r<b(z)\}$.
Then
\begin{equation}\label{eq:wall-in-U}
 \Gamma_a:=\{(a(z),z):z\in J\}\subset\mathcal U,
 \qquad \Gamma_a\text{ has a neighborhood in }\mathcal U.
\end{equation}
Write
\[
 E=\{t\in J:(0,t)\in\overline{Z_N\cup Z_S}\}
\]
for the axis contact heights. Proposition~\ref{prop:finite-poles}
gives an accessible nodal arc in $D_a$ at every $t\in E$.

\smallskip\noindent
\textit{Step 2: a cap crosscut clears an axis tail.}
Suppose a component $\mathcal C$ of $(Z_N\cup Z_S)\cap D_b$
has an axis end $(0,t)$, $t\in J$, and meets $\mathcal U$.
Then $\mathcal C\subset Z_N\cap\mathcal U$.
Take disjoint short end-arcs $\eta_0,\eta_1$ entering $\mathcal U$
from $(0,t)$ and $p$, respectively, using the nodal arc and
the wall $\Gamma_a$. Join their interior endpoints by a simple
polygonal arc $Q\subset\mathcal U$. With
\[
 s_i=\min\{s:\eta_i(s)\in Q\}>0,
\]
the two prefixes and the intervening subarc of $Q$ form an embedded
arc $\gamma$ from $(0,t)$ to $p$ with
$u(\gamma\setminus\{(0,t),p\})\subset B_{\alpha_0}(N)$.
Let $V$ be the Jordan domain bounded by $\gamma$ and
$\{0\}\times[t,1]$. Since $D_b$ is a Jordan domain, $V\subset D_b$.

Apply the cap fold $\Pi_{N,\alpha}$ from
Subsection~\ref{sec:cap-folds} in $\widehat V$.
It fixes the interface and gives
\[
 \begin{gathered}
 E(v)\le E(u),\qquad q_v=1\text{ on }(t,1),\qquad
 q_v=q\text{ on }(-1,t),\\
 L(v)-L(u)=-k\int_t^1\mathbf1_{\{q=-1\}}\,dz\le0.
 \end{gathered}
\]
Equation~\eqref{eq:boundary-equality} gives $u=v$ off the axis.
Every point $(0,s)$ with $t<s<1$ has a meridian neighborhood
inside $V$, so Lemma~\ref{lem:cap-removal} yields
\[
 \{(0,0,s):t<s<1\}\subset\Reg(u).
\]
The assertion follows in this case.

\smallskip\noindent
\textit{Step 3: otherwise every contact uses a distinct bottom incidence.}
Assume the alternative in Step~2 does not occur.
Every component accessible from $E$ is then disjoint from
$\mathcal U$ and cannot cross $\Gamma_a$.
Continue its accessible arc without backtracking, as in the proof
of Proposition~\ref{prop:finite-poles}. Until its first exit
through $\Gamma_0$, it remains in $D_a$.
If no exit occurs, local nodal finiteness, including the finite
stars at $\mathcal E_0$, gives a connected tail-limit set
\[
 \varnothing\ne K_\infty
 =\bigcap_{T>0}\overline{\gamma([T,\infty))}
 \subset(\{0\}\times E)\cup\{p\}.
\]
The good height $z_0$ excludes the lower axis endpoint, and
$E\cap\mathcal G=\varnothing$ excludes nontrivial axis intervals.
Hence $K_\infty$ is a singleton. The completed path has two
closed-axis ends and is forbidden by \eqref{eq:boundary-half-turn}.
Thus every contact reaches a bottom half-branch.

The same forest argument as in Proposition~\ref{prop:finite-poles}
shows that two distinct contacts cannot use the same half-branch:
their paths would form a forbidden arc with two axis ends.
Consequently,
\[
 E\hookrightarrow\{\text{nodal half-branches incident to }\mathcal E_0\},
 \qquad e:=\#E\le M.
\]
On each component of $J\setminus E$, the compact-interval
argument of Lemma~\ref{lem:interior} applies without change.
Proposition~\ref{prop:two-points} therefore gives
\[
 \#\{z\in J:(0,0,z)\in\Sing(u)\}
 \le e+2(e+1)\le3M+2<\infty.\qedhere
\]
\end{proof}

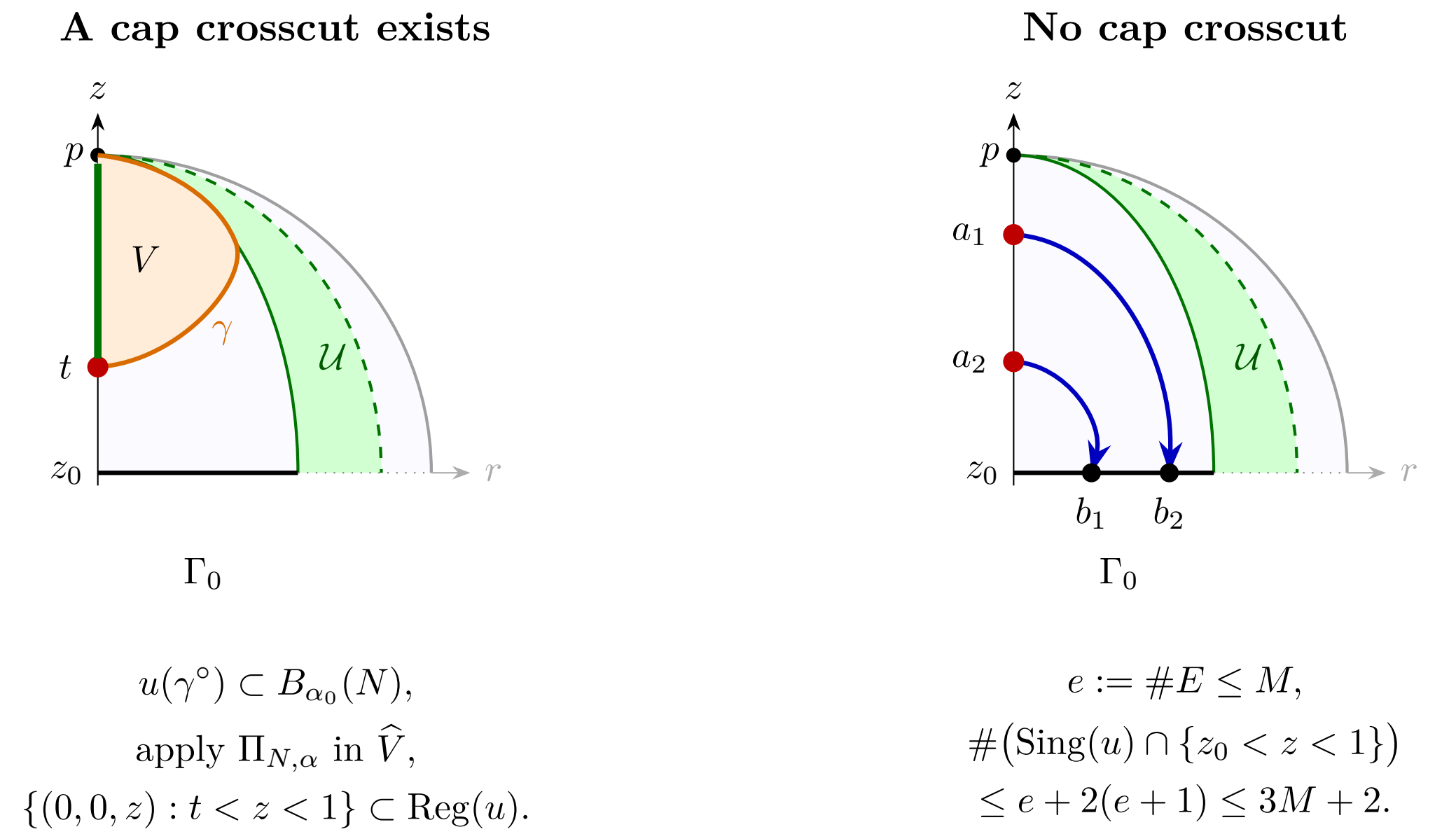
\begin{figure}[htbp]
\centering
\tikzset{
  pics/boundarycollar/.style={code={
    \fill[blue!2]
      (0,0)--(0,3)
      arc[start angle=90,end angle=0,radius=3]--cycle;

    \fill[green!18]
      (0,3)
      arc[start angle=90,end angle=0,
          x radius=2.55,y radius=3]
      --(1.8,0)
      arc[start angle=0,end angle=90,
          x radius=1.8,y radius=3]
      --cycle;

    \draw[gray!75,thick]
      (0,3) arc[start angle=90,end angle=0,radius=3];
    \draw[green!45!black,thick,dashed]
      (0,3)
      arc[start angle=90,end angle=0,
          x radius=2.55,y radius=3];
    \draw[green!45!black,thick]
      (0,3)
      arc[start angle=90,end angle=0,
          x radius=1.8,y radius=3];

    \draw[->] (0,-.12)--(0,3.40) node[above] {$z$};
    \draw[very thick] (0,0)--(1.8,0);
    \draw[gray,dotted] (1.8,0)--(3,0);
    \draw[->,gray!65] (3,0)--(3.35,0) node[right] {$r$};

    \fill (0,3) circle (2pt);
    \node[left] at (0,3) {$p$};
    \node[left] at (0,0) {$z_0$};
    \node at (.95,-.95) {$\Gamma_0$};
    \node[green!35!black] at (2.12,1.1) {$\mathcal U$};
  }},
  axiscontact/.style={
    circle,fill=red!75!black,inner sep=2pt
  },
  bottomport/.style={
    circle,fill=black,inner sep=1.8pt
  }
}

\begin{minipage}[t]{.48\textwidth}
\centering
\textbf{A cap crosscut exists}

\medskip
\begin{tikzpicture}[x=1.05cm,y=1cm,font=\small,>=Stealth]
  \pic {boundarycollar};

  \fill[orange!15]
    (0,1)
    .. controls (.75,1.05) and (1.40,1.85) .. (1.23,2.20)
    .. controls (1.05,2.65) and (.50,2.95) .. (0,3)
    --cycle;

  \draw[orange!85!black,very thick]
    (0,1)
    .. controls (.75,1.05) and (1.40,1.85) .. (1.23,2.20)
    .. controls (1.05,2.65) and (.50,2.95) .. (0,3);

  \node[orange!85!black] at (1.12,1.32) {$\gamma$};
  \node at (.43,2.02) {$V$};
  \node[axiscontact,label=left:{$t$}] at (0,1) {};
  \draw[green!45!black,line width=2pt]
    (0,1.08)--(0,2.92);
\end{tikzpicture}

\small
\[
 \begin{gathered}
 u(\gamma^\circ)\subset B_{\alpha_0}(N),\\
 \text{apply }\Pi_{N,\alpha}\text{ in }\widehat V,\\
 \{(0,0,z):t<z<1\}\subset\Reg(u).
 \end{gathered}
\]
\end{minipage}
\hfill
\begin{minipage}[t]{.48\textwidth}
\centering
\textbf{No cap crosscut}

\medskip
\begin{tikzpicture}[x=1.05cm,y=1cm,font=\small,>=Stealth]
  \pic {boundarycollar};

  \draw[blue!75!black,very thick,->]
    (0,2.25)
    .. controls (.85,2.20) and (1.45,1.10) .. (1.40,0);
  \draw[blue!75!black,very thick,->]
    (0,1.05)
    .. controls (.48,1.02) and (.82,.48) .. (.70,0);

  \node[axiscontact,label=left:{$a_1$}] at (0,2.25) {};
  \node[axiscontact,label=left:{$a_2$}] at (0,1.05) {};
  \node[bottomport,label=below:{$b_1$}] at (.70,0) {};
  \node[bottomport,label=below:{$b_2$}] at (1.40,0) {};
\end{tikzpicture}

\small
\[
 \begin{gathered}
 e:=\#E\le M,\\
 \#\bigl(\Sing(u)\cap\{z_0<z<1\}\bigr)\\
 \le e+2(e+1)\le3M+2.
 \end{gathered}
\]
\end{minipage}

\caption{The boundary-pole alternative.
The green strip between $r=a(z)$ and $r=b(z)$ lies in the cap
component $\mathcal U$, which may extend beyond that strip.
Left: $\gamma$ and the axis segment bound $V$; the cap comparison
makes the green axis tail regular.
Right: accessible pole-level paths must reach the artificial
section $\Gamma_0$, with distinct contacts using distinct bottom
half-branches. The number $M$ counts incidences, including tangencies.
Red points are axis contacts, not necessarily singularities.}
\label{fig:boundary-pole-alternative}
\end{figure}

\begin{proof}[Proof of Theorem~\ref{thm:global}]
Existence follows from Lemma~\ref{lem:existence}.
For any minimizer, Lemmas~\ref{lem:canonical-lift}
and~\ref{lem:interior} give interior local finiteness and
$\Sing(u)\subset\axis$. After dilation to $R=1$, apply
Proposition~\ref{prop:boundary-nonaccumulation} at $p_+$ and,
by $z\mapsto-z$, at $p_-$.
Choose $-1<t_-<t_+<1$ so that both outer axis tails contain
only finitely many singularities. The remaining segment
\[
 K=\{(0,0,z):t_-\le z\le t_+\}\Subset B_1
\]
has finitely many singularities by Lemma~\ref{lem:interior}.
Thus
\[
 \#\Sing(u)
 =\#(\Sing(u)\cap K)
 +\#(\Sing(u)\cap\{z<t_-\})
 +\#(\Sing(u)\cap\{z>t_+\})<\infty.
\]
\end{proof}

\section*{Acknowledgments}

This work is partially supported by the National Key R\&D Program of China
under Grant 2023YFA1008801 and NSF of China under Grant No. 12271476.
\subsection*{The use of AI tools}
The authors acknowledge the assistance of AI tools (ChatGPT) in
preparing this article. ChatGPT 6.0 Astra helped identify a suitable
proof strategy, while ChatGPT 6.1 Sol assisted in simplifying the
proofs, improving the exposition connecting the interior and
boundary arguments, and producing four color schematic figures
illustrating the proof strategy. All results and proofs have been
carefully checked by the authors, who take full responsibility
for the content of the article.

\bibliographystyle{amsplain}
\bibliography{hlp_extension}
\end{document}